\documentclass[10pt]{amsart}
\usepackage{graphicx,url}
\usepackage{indentfirst}
\usepackage{amsfonts}
\usepackage{amsmath}
\usepackage{amssymb}
\usepackage{color}

\newcommand{\al} {\alpha}

\newtheorem{obs}{Remark}[section]

\newtheorem{prop}[obs]{Proposition}
\newtheorem{teo}[obs]{Theorem}
\newtheorem{lem}[obs]{Lemma}

\def\dem{{\bf Proof:}\hspace{.2in}}

\def\dem{{\bf Proof:}\hspace{.2in}}

\def\demT{{\bf Proof:}\hspace{.2in}}

\def\demP215{{\bf Proof of Proposition 1.6:}\hspace{.2in}}

\begin{document}
\author{
J. Bastos,  A. Messaoudi, D. Smania , T. Rodrigues}

\address{Departamento de Matemática, UNESP - Universidade Estadual Paulista, Rua Crist\'ov\~ao Colombo, 2265, Jardim Nazareth, 15054-000, S\~ao Jos\'e do Rio Preto, SP, Brazil.}
\email{messaoud@ibilce.unesp.br,\ jeferson@ibilce.unesp.br}

\address{Departamento de Matem\'atica, ICMC-USP, Avenida do Trabalhador S\~ao Carlense, 400, Caixa Postal 668, 13560-970, S\~ao Carlos, SP, Brazil.}
\email{smania@icmc.usp.br}

\address{Departamento de Matem\'atica, UNESP - Universidade Estadual Paulista, AV. Eng. Luiz Ed. Carrijo Coube, 14-01, Vargem Limpa, 17033-360, Bauru, SP, Brazil.}
\email{tatimi@fc.unesp.br}

\date{\today}
\title{A class of cubic Rauzy Fractals}
\begin{abstract} In this paper, we study arithmetical and topological properties for  a class of
Rauzy fractals ${\mathcal R}_a$ given by the polynomial $x^3- ax^2+x-1$ where $a \geq
2$ is an integer.
In particular, we prove the number of neighbors of ${\mathcal R}_a$ in the periodic tiling is  equal to $8$. We also
 give explicitly an automaton that generates the boundary of ${\mathcal R}_a$. As a consequence, we prove that ${\mathcal R}_2$ is homeomorphic to a topological disk. \end{abstract}
\subjclass[2000]{}

\thanks{A.M. was supported by Brasilian CNPq grant
305939/2009-2}

\thanks{D.S. was partially supported by CNPq  303669/2009-8 305537/2012-1 and FAPESP  2008/02841-4, 2010/08654-1 }


\maketitle

\section{Introduction}

In 1982, G. Rauzy \cite{rauzy1} defined a compact subset ${\mathcal E}$
of $\mathbb{C}$ called classical Rauzy Fractal as
$${\mathcal E}= \{\sum_{i=0}^{+\infty} \varepsilon_i \alpha^i,\; \varepsilon_i \in \{0,1\},\varepsilon_i \varepsilon_{i+1} \varepsilon_{i+2} \ne 111,
\ \forall i\geq 0\},$$ where $\alpha$ is one of the two complex roots
of modulus $<1$ of the polynomial $P(x)= x^3- x^2-x-1$.

The classical Rauzy fractal has many beautiful properties: It is a
connected set, with interior simply connected, and boundary fractal.
Moreover, it induces a periodic tiling of the plane $\mathbb{C}$
modulo the group $\mathbb{Z} \alpha^{-3}+ \mathbb{Z} \alpha^{-2}.$

The Rauzy fractal was studied by many mathematicians  and was
connected to  to many topics as: numeration systems
(\cite{ali2},\cite{ali22}, \cite{P99}), geometrical representation
of symbolic dynamical system (\cite{arnoux/ito}, \cite{arnoux/sano}, \cite{arnoux/rauzy},
\cite{canterini}, \cite{HZ98}, \cite{ali1}, \cite{pytheas}, \cite{T06},  \cite{S96}),
multidimensional continued fractions and simultaneous approximations
(\cite{ABI02}, \cite{CHM01}, \cite{C99}, \cite{HM06}), auto-similar
tilings (\cite{akiyama1}, \cite{akiyama}, \cite{arnoux/ito},
\cite{P99}) and Markov partitions of Hyperbolic automorphisms of
Torus (\cite{KV98}, \cite{ali1}, \cite{P99}).

There are many ways of constructing Rauzy fractals, one of them is by
$\beta$-expansions.

Let $\beta > 1$ be a real number and $x\in\mathbb{R}^{+}$. Using
greedy algorithm, we can  write $x$ in base $\beta$ as  $x=
\sum_{i=-\infty}^{k}a_i \beta^{i}$ where $  k \in \mathbb{Z}$  and
$a_{i}$ belong to the set $A$ where  $A=\{0, \ldots, \beta-1\}$ if
$\beta \in \mathbb{N}$ or $A=\{0, \ldots, \lfloor \beta \rfloor\}$
otherwise, where $\lfloor \beta \rfloor$ is the integer part of
$\beta$. The sequence $(a_i)_{i \leq k}$ is called $\beta-$
expansion of $x$ and is also denoted by $a_k a_{k-1}\ldots$ The
greedy algorithm can be defined as follows (see \cite{parry} and
\cite{frougny1}): denote by $\{y\}$ the fractional party of a number
$y$. There exists an integer
 $k \in
\mathbb{Z}$ such  $\beta^{k} \leq x <\beta^{k+1}$. Let $x_k=
\lfloor x/\beta ^{k}\rfloor$ and  $r_k = \{ x/\beta ^{k}\}$. Then
for $i < k,$ put $x_i = \lfloor \beta r_{i+1} \rfloor$ and $r_i =
\{\beta r_{i+1}\}$. We get $$x= x_{k} \beta^{k}+ x_{k-1} \beta^{k-1}+
\cdots$$ if $k <0 \; (x <1)$, we put $x_{0}=x_{-1}= \cdots =
x_{k+1}=0.$ If an expansion $(x_{i})_{i \leq k}$ satisfies $x_i=0$ for all $i <n$, it is said to be finite and the ending zeros are omitted.
It will be denoted by $(x_{i})_{ n \leq i \leq k}$ or $x_k \ldots x_n$.

 Now, assume that $\beta$ is a Pisot number of degree $d \geq 3$, that means that $\beta$ is an algebraic integer of degree $d$ whose Galois' conjugates have modulus less than one.
 We denote by $b_2, \ldots, \beta_r$ the real Galois conjugates of $\beta$ and by $\beta_{r+1},  \ldots, \beta_{r+s}, \beta_{r+s+1} =\overline{\beta_{r+1}},\ldots, \beta_{r+2s}= \overline{\beta_{r+s}}$ its complex Galois conjugates.
 Let $\psi= (\beta_2,\ldots, \beta_{r+s}) \in \mathbb{R}^{r-1}\times \mathbb{C}^s$ and put $\psi^{i}= (\beta_2^{i},\ldots, \beta_{r+s}^{i}) $ for all $ i \in \mathbb{Z}.$

The Rauzy fractal is by definition the set

 $$\mathcal{R}={\mathcal R}_{\beta}= \{\sum_{i=0}^{+ \infty} a_i \psi^{i},\; (a_{i})_{i \geq 0 } \in E_{\beta}\},$$
 where
 $$E_{\beta}=
\{ (x_{i})_{i \geq k },\; k \in \mathbb{Z} \; \vert \;
 \forall n \geq k,\; (x_{i})_{n \geq i \geq k} \mbox { is a finite
 } \beta \mbox { expansion } \}.$$

Observe that  ${\mathcal R}_{\beta}$ is a compact subset of $\mathbb{R}^{r-1} \times \mathbb{C}^s \approx \mathbb{R}^{d-1}$.

 For example, if $\beta >1$ is a root  of the polynomial $P(x)= x^3- x^2-x-1$, we obtain the classical Rauzy fractal ${\mathcal R}_{\beta}={\mathcal E}$.

 An important class of Pisot numbers are those such that the associated Rauzy fractal has $0$ as an interior point.
This  numbers were characterized by Akiyama in \cite{akiyama2}. They
are exactly the Pisot numbers
 that satisfy
$$\mathbb{Z}[\beta]\cap [0, +\infty[ \subset  \mbox{Fin}(\beta) \mbox { (called property (F)) },$$
where $\mbox{Fin}(\beta)$ is the set of nonnegative real numbers
which have a finite $\beta$-expansion.

In this paper we study properties of the Rauzy fractal associated to
a class of cubic unit Pisot numbers that satisfy property (F). These
numbers were characterized in \cite{akiyama} as being exactly the
set of dominant roots of the polynomial (with integers coefficients)
$$P_{a,b}(x)=x^3-ax^2-bx-1,\ a\geq 0,\ -1\leq b \leq a+1.$$
(If $b=-1$ add the restriction $a\geq 2$).

%
%

%
%
%
%
%

%
%
%
%
%
%
%

In particular, this set divided into three subsets:
\begin{itemize}
\item[a)]  $0 \geq b \geq a$, and in this case $d(1,\beta)
= \cdot a b 1.$
\item[b)] $b=-1,\;
a \geq 2.$ In this case   $d(1,\beta) = \cdot (a-1)(a-1)01.$
\item[c)]
$b=a+1,$  and in this case $d(1,\beta) = \cdot (a+1) 0 0 a 1$, where
$d(1, \beta)$ is the R\'enyi $\beta$-representation of $1$ (see
\cite{rényi}).
\end{itemize}

Geometrical and arithmetical properties of the Rauzy fractal
associated to polynomials $P_{a,b},\ a\geq b\geq 1$ were studied in
\cite{TLMS}. Here we will study the case $a \geq2, \ b=-1$. In this case  the polynomial $p(x)= x^3 -ax^2+x-1= (x -\beta)(x-\alpha) (x-
\gamma)$, where $\beta >1$ and
 $ \alpha, \gamma \in \mathbb{C} \setminus \mathbb{R}$,
 and  the
Rauzy fractal
$$
\mathcal{R}_a=\left\{\sum_{i=0}^{\infty} a_{i}
\alpha^{i},\; a_{i}a_{i-1}a_{i-2}a_{i-3}<_{lex} d(1,
\beta)=(a-1)(a-1)01 ,\; \forall i \geq 0 \mbox{ where } a_{-1}=a_{-2}= a_{-3}=0 \right\},$$ where $<_{lex}$
is the lexicographic order on finite words.

On the other hand, consider the sequence  $ R_0  = 1, \ R_1 = a, \; R_2 = a^2,\; R_{n+3}
= a  R_{n+2} - R_{n+1} + R_{n} \; \forall n \geq 0.$
It is known, using greedy algorithm that for all nonnegative integer $n$ can can be written as $n= \sum_{i=0}^{N} a_{i}R_i$.
The sequence $(a_i)_{0\leq i\leq N}$ is called a greedy R-expansion.

The Rauzy fractal is equal
$$
{\mathcal R}_a=\left\{\sum_{i=0}^{\infty} a_{i} \alpha^{i},\ \forall N\geq
0\ (a_i)_{0\leq i\leq N}\ \mbox{is a greedy R-expansion}\right\}.$$

We will also study properties of another set very closed to the Rauzy Fractal. We call this set
the $G$-Rauzy fractal and define it by

$$
\mathcal{G}_a=\{\sum_{i=0}^{\infty} a_{i} \alpha^{i}, \forall N \geq 0
,\ (a_{i})_{0 \leq i \leq N} \mbox{ is a greedy $ G$-expansion
}\},$$ where $G= (G_n)_{n \geq 0}$ where $  G_0  = 1, \ G_1 = a, \
G_2 = a^2 +b,\; G_{n+3} = a  G_{n+2} +  b G_{n+1} + G_{n} \; \forall
n \geq 0.$

The set $\mathcal{G}$ was defined
in \cite{HM06} by Hubert and Messaoudi. They used it to prove that $(G_n)_{n \geq 0}$ is the sequence of  best approximations of the vector  $(1 /\beta,\; 1 / \beta^{2})$  (for a certain norm on $\mathbb{R}^2$ called the Rauzy norm ${\mathcal
N}).$

In the case where $b=-1$ and $a \geq 2$ it is known (see
\cite{HM06}) that the set of $ G$-expansions is equal to the set of
$(\varepsilon_{i})_{0 \leq i \leq N}$ that satisfy the following
conditions:
$$
\varepsilon_{i}\varepsilon_{i-1}\varepsilon_{i-2}\varepsilon_{i-3} <_{lex} d(1, \beta)= (a-1)(a-1)01,\ \forall i \geq 3,$$
 and the initial conditions
 $$
 \varepsilon_{0}< a,\ \varepsilon_{1}\varepsilon_{0}<_{lex} (a-1)(a-1),\ \varepsilon_{2}\varepsilon_{1}
\varepsilon_{0}<_{lex} (a-1)(a-1)0.$$

Observe that the above initial conditions from the fact that:
$\varepsilon_{0}G_0 < G_1=a,\; \varepsilon_{0}G_0+ \varepsilon_{1}G_1 <_{lex} G_2= (a-1)G_1+ (a-1)G_0,\;
\varepsilon_{0}G_0+ \varepsilon_{1}G_1 +\varepsilon_{2}G_2 < G_3= (a-1)G_2+ (a-1)G_1.$


%
%

%

Many topological properties of ${\mathcal R}_a$ are known (see
\cite{akiyama,IK91,ali1,ali22,rauzy1}): It's a connected compact
subset of $\mathbb{C}$, with interior simply connected and fractal
boundary, moreover it  induces a periodic tiling of the plane modulo
$\mathbb{C}$. It can be also seen as geometrical realization of
the dynamical system associated to the substitution $\sigma$ defined
by: $\sigma(1)= 1^{a-1}2,\; \sigma(2)= 1^{a-1}3,\; \sigma(3)=4,\;
\sigma(4)=1$.

To our knowledge,  geometrical and topological properties of the set
${\mathcal G}_a $ were not yet studied. In this paper, we show that
$\mathcal{G}_a$  induces a periodic tiling of the complex plane.  We
also construct an explicit finite state automaton ${\mathcal A}$ that
generates both boundaries of ${\mathcal R}_a$ and ${\mathcal G}_a$. With
this we prove that for all $a \geq 2,\; {\mathcal R}_a$ has $8$
neighbors while ${\mathcal G}_a $ has $6$ neighbors (in the periodic
tiling). The interest of giving explicitly the automaton
$\mathcal{A}$ remains in the fact that the study of properties of
$\mathcal{A}$
 give topological and metrical  information about the
boundary $\mathcal{G}_{a}$ and  $\mathcal{R}_{a}$.

Here, we prove that the boundary of $\mathcal{R}_{2}$ is homeomorphic to a topological circle. This study can be done for all integer $a \geq 2$.

The paper is divided by the following manner. In the second section, we give some notations. In the third section, we study some properties of the boundary of ${\mathcal G}_a$,
in the fourth section, we construct an explicit finite state automaton that recognizes the boundaries of $ \mathcal{G}_a$ and $\mathcal{R}_a$ for all $a \geq 2$.
The fifth section is devoted to the study topological properties of the boundary of $\mathcal{R}_2$. In particular, using the automaton, we prove that the boundary of $\mathcal{R}_2$ is homeomorphic to a circle.

\section{Notations and definitions}
Denote by $E(G)$ (resp. $E(R)$) the set of sequences  $(a_n)_{n \in
\mathbb{Z}}$ belonging to  $\{0,1,\ldots, a-1\}^{\mathbb{Z}}$ such
that, there exists an integer $k \in \mathbb{Z}$ satisfying $a_k >0$ and $a_n=0$
for all  $n < k$, moreover  for all $p \geq k$, the sequence $(a_n)_{k
\leq n \leq p}$ is a $G$-expansion (resp. $R$-expansion ).
That is

$E(R)= \{(a_{n})_{n \in
\mathbb{Z}}, \; \exists k \in \mathbb{Z},\; a_k >0,\;  a_i = 0 \mbox { for all } i <k,\;  a_{i}a_{i-1}a_{i-2} a_{i-3} <_{lex} (a-1)(a-1)01,\; \forall i \geq k  \} ,$
and
$E(G)= \{(a_{n})_{n \in
\mathbb{Z}}, \; \exists k \in \mathbb{Z},\; a_k >0, \;a_i = 0 \mbox { for all } i <k,\;  a_{i}a_{i-1}a_{i-2} a_{i-3} <_{lex} (a-1)(a-1)01,\; \forall i \geq k,\;
 a_{k}< a,\, a_{k+1}a_{k}<_{lex} (a-1)(a-1),\ a_{k+2}a_{k+1}
a_{k}<_{lex} (a-1)(a-1)0  \} .$
Observe
that $E(G) \subset E(R).$

\noindent  We will identify a sequence   $(a_n)_{n \in \mathbb{Z}}$
belonging to   $E(R)$  such that $a_n= 0$ for all $ n <
k$ with the sequence $(a_n)_{n \geq k}.$

Let $(a_n)_{n \geq k}$ be an element of  $E(R)$. Assume
that there  exists $p \in \mathbb{Z}$ such that  for all $n
> p,\; a_n= 0.$ This sequence will be denoted by $(a_n)_{k \leq n \leq
p}$.

For technical reasons, we will consider
$$\mathcal{R}_a=\left\{\sum_{i=2}^{\infty}a_i\al^{i}, a_i a_{i-1}a_{i-2}a_{i-3}<_{lex}(a-1)(a-1)01, \forall i\geq
2, \mbox{ where } a_{1}=a_{0}= a_{-1}=0\right\} $$ and
$$\mathcal{G}_a=\left\{\sum_{i=2}^{\infty}a_i\al^{i},\ \forall\ N\geq 2\  ,(a_i)_{2\leq i \leq N}\ \mbox{is a Greedy G-expansion}\right\}.$$
\section{Properties of $\mathcal{G}_a$ and its boundary}

\begin{teo}
\label{R}
 The set $\mathcal{G}_a$ induces a periodic tiling of the complex plane, that is,
\begin{enumerate}
\item[a)] $\mathbb{C}=\bigcup_{u \in
\mathbb{Z}+\mathbb{Z}\alpha}(\mathcal{G}_a+u)$;
\item[b)] $int(\mathcal{G}_a+u) \cap (\mathcal{G}_a+v) \neq \emptyset,\ u,\
v \in \mathbb{Z}+\mathbb{Z}\alpha$ implies que  $u=v$.
\end{enumerate}\label{teo23}
\end{teo}

\begin{obs}
 The proof can be deduced from \cite{rauzy1} (done in case of  Rauzy fractal ${\mathcal E}$,
see also \cite{DF}). For clarity, we will give the proof here.
\end{obs}

\hfill $\Box$

Consider  the sequence $G'= (G'_{n})_{n \geq 0}$ by $G'_{0}=0,\
G'_{1}=0,\ G'_{2}=1,\ G'_{n+3}=a G'_{n+2}-G'_{n+1}+G'_{n},\forall\ n
\geq 0$. Then $G_{n}= G'_{n+2}$ for all integer $n \in \mathbb{N}$.


\begin{prop}\label{prop12} The following properties are valid:
\begin{itemize}
\item[i)] All natural integer $n$ can be written by unique way as $n=\sum_{i=2}^{N} \varepsilon_{i} G'_{i}$ where
$(\varepsilon_{i})_{2 \leq i \leq N} \in E(G)$.

\item [ii)] Let $(a_{i})_{l \leq i \leq N}$  and $(b_{i})_{l^{'} \leq i \leq \infty}$ be two elements of  $E(R)$ (resp. $E(G)$)
 such that $a_{l} >0$ and $b_{l^{'}}>0.$ If
$\sum_{i=l}^{N}a_{i}\alpha^{i} =
\sum_{i=l^{'}}^{\infty}b_{i}\alpha^{i}$ then $l=l^{'}$ and for all
$i \geq N,\ b_{i}=0$ and for all $l \leq i \leq N,\ a_{i}=b_{i}.$

\item[iii)] Let $(\varepsilon_{i})_{2 \leq i \leq N} \in E(R)$ (resp. $  E(G)$ ) then $\sum_{i=2}^{N} \varepsilon_{i}\alpha^{i}
\in int(\mathcal{R})$ (resp. $int(\mathcal{G})$). In particular, $0 \in int(\mathcal{R})$ (resp. $0 \in int(\mathcal{G})$).

\item[iv)] Let $z \in \mathbb{Z}[\beta]\cap\mathbb{R}^{+}$ then there exist a sequence
$(a_{i})_{k \leq i \leq l} \in E(R),\ k \leq l$ such that
$z=\sum_{i=k}^{l}a_{i}\beta^{i}$.

\item[v)] For all $n \geq 2$ we have $\beta^{n}=G'_{n}\beta^{2}+(G'_{n-2}-G'_{n-1})\beta + G'_{n-1}$.
In particular if $(\varepsilon_{i})_{2 \leq i \leq l} \in E(G)$ then
$\sum_{i=2}^{l} \varepsilon_{i}\beta^{i}= n\beta^{2}+r(n)\beta+s(n)$
where $n=\sum_{i=2}^{l} \varepsilon_{i}G'_{i}, r(n)= \sum_{i=2}^{l}
\varepsilon_{i}(G'_{i-2}-G'_{i-1})$ and $s(n)=\sum_{i=2}^{l}
\varepsilon_{i}G'_{i-1}.$

\item[vi)] Let $(a_{i})_{l \leq i \leq k}$ and $(b_{i})_{l \leq i \leq k}$ be elements of $ E(G)$ (resp.$E(R)$ ). Then $\sum_{i=l}^{k} a_{i} \beta^{i} < \sum_{i=l}^{k} b_{i} \beta^{i}$ if, only if
        $(a_{i})_{l \leq i \leq k} <_{lex} (b_{i})_{l \leq i \leq k}$.

\item[vii)] Let $c,\ d \in \mathbb{R}$ such that  $\alpha^{2}=c+d\alpha,$ then $1,\ c$ and $d$ are  $\mathbb{Q}$-Linearly independent.

\end{itemize}
\end{prop}

\begin{obs}
The results given in Proposition \ref{prop12} are classical. For i) and vi), see \cite{L}. For ii) see \cite{HM06}. The results iii) and vii) can be found in \cite{akiyama}.

For iv), see \cite{SF}.  v) is left to the reader and can bem done by induction.

\end{obs}

{\bf Proof of Theorem 2.1.}\\

 Let $z \in \mathbb{C}$ and $\epsilon > 0$.
Using item (vii) of Proposition \ref{prop12} and Kronecker's
Theorem, we deduce that the set $\{n \alpha^{2}+p\alpha +q,\ n \in
\mathbb{N},\ p,\ q \in \mathbb{Z}\}$ is  dense in $\mathbb{C}$.
Then there exists a sequence $(z_{k})_{k \geq 0} \in \mathbb{C}$
such that
\[
z_{k}=n_{k} \alpha^{2} + p_{k}\alpha+r_{k},\ n_{k} \in \mathbb{N},\
p_{k},\ q_{k} \in \mathbb{Z}
\]
and for all $k \geq k_{0},\ \mid z_{k} - z \mid < \epsilon$. Let
$A_{k}= n_{k} \alpha^{2} + r(n_{k}) \alpha + s(n_{k}),$ where
$r(n_{k})$ and $s(n_{k})$ are defined in item (v) in Proposition
\ref{prop12}. We have $A_{k} \in \mathcal{G}_a$.

On the other hand,
\[
A_{k}=n_{k} \alpha^{2} + p_{k}\alpha+r_{k} + (r(n_{k})-p_{k})\alpha+
(s(n_{k})-r_{k})=z_{k} + t_{k}\alpha+m_{k},
\] where
$t_{k}=r(n_{k})-p_{k}$ e $m_{k}=s(n_{k})-r_{k}$. Then, $z_{k} +
t_{k}\alpha+m_{k} \in \mathcal{G}_a$.

On the other hand, for all $k \geq k_{0},$
\[
\mid z+t_{k}\alpha +m_{k} \mid \leq \mid z - z_{k} \mid + \mid z_{k}
+ t_{k} \alpha + m_{k} \mid < \epsilon + d,
\]
where $d$ is the diameter of $ \mathcal{G}_a$.
 Since $\mathbb{Z}+\mathbb{Z}\alpha$ is a lattice,
 there exists an increasing sequence
$(k_{i})_{i \geq 1}$ of integer numbers such that for all $i,\ j \in
\mathbb{N},$ \ $t_{k_{i}} \alpha+r_{k_{i}} = t_{k_{j}} \alpha +
r_{k_{j}}.$ Then, there exist $t,\ r \in \mathbb{Z}$ such that
$t_{k_{i}}=t_{k_{j}}=t$ and $r_{k_{i}}=r_{k_{j}}=r,$ for all $i,j
\in \mathbb{N}$. As $z_{k_{i}}+t_{k_{i}}\alpha+m_{k_{i}} \in
\mathcal{G}_a,\  \lim_{i \to +\infty} z_{k_{i}} = z$ and $\mathcal{G}_a$ is a
closed set, we have that $z+t\alpha+r \in \mathcal{G}_a.$   \hfill$\Box$

\vspace{0.3cm}

To prove, item b),  it is sufficient to establish that if
$(int(\mathcal{G}_a)+u) \cap \mathcal{G}_a \neq \emptyset$ where $u \in
\mathbb{Z}+\mathbb{Z}\alpha$ then $u=0$.

Assume  that there exist $p,\ q \in \mathbb{Z}$ and an
element $z=\sum_{i=2}^{\infty} \varepsilon_{i}\alpha^{i} \in
\mathcal{G}_a$ such that $z+p+q\alpha \in int(\mathcal{G}_a)$. Thus
there is an integer $n_{0} \geq 0$ such that for all $n \geq n_{0}$
\begin{equation}\label{eq001}
\sum_{i=2}^{n} \varepsilon_{i}\alpha^{i} + p +q\alpha \in
\mathcal{G}_a.
\end{equation}
\textbf{Case 1:} The set $\{i \geq 2, \varepsilon_{i} \neq 0\}$ is
infinite.

In this case, as $\beta > 1,$ then there exists a integer $N \geq
n_{0}$ such that $\sum_{i=2}^{N} \varepsilon_{i} \beta^{i} +
p+q\beta > 0.$ By item (iv) of Proposition \ref{prop12} we deduce that
\begin{equation}\label{eq002}
\sum_{i=2}^{N} \varepsilon_{i}\beta^{i} + p+q\beta = \sum_{i=l}^{M}
d_{i}\beta^{i},\ \mbox{where} \hspace{0.2cm} (d_{i})_{l \leq i \leq
M} \in  E(R),\ l,\ M \in \mathbb{Z}.
\end{equation}

From (\ref{eq001}) and (\ref{eq002}) we have that $\sum_{i=l}^{M}
d_{i}\alpha^{i} = \sum_{i=2}^{\infty} e_{i}\alpha^{i} \in
\mathcal{G}_a$.

Therefore, from item  (ii) of  Proposition \ref{prop12}, we have
$e_{i}=0$ for all $i > M.$ Then, \[
\begin{array}{ccc}
\sum_{i=2}^{N} \varepsilon_{i}\beta^{i}+p + q \beta &=&
\sum_{i=2}^{M}e_{i}\beta^{i}\\ \ \\&=& \sum_{i=l}^{M} d_{i}
\beta^{i}.
\end{array}
\]

According to item (v) of Proposition \ref{prop12}, we have
\[
\widetilde{n}
\beta^{2}+(r(\widetilde{n})+q)\beta+(s(\widetilde{n})+p) =
\widetilde{l}\beta^{2}+r(\widetilde{l})\beta+s(\widetilde{l}).
\]
where $\widetilde{n}=\sum_{i=2}^{N}\varepsilon_{i}G'_{i}$ and
$\widetilde{l}=\sum_{i=2}^{M}e_{i}G'_{i}$. Therefore,
$\widetilde{l}=\widetilde{n}$ and $\varepsilon_{i}=e_{i}$ for all
$i$ (by (i) of Proposition \ref{prop12}). Thus, $p=q=0$.

\textbf{Case 2:} The set $\{ i \geq 2,\ \varepsilon_{i} \neq 0 \}$
is finite.

Let $N=max\{i \geq 2, \varepsilon_{i} \neq 0\}$. If
$\sum_{i=2}^{N}\varepsilon_{i}\beta^{i} + p + q\beta \geq 0$ then we
use the same argument than in case 1.

Assume that $\sum_{i=2}^{N}\varepsilon_{i}\beta^{i} + p +q\beta <
0.$ We have
\[
\sum_{i=2}^{N} \varepsilon_{i}\alpha^{i} + p + q\alpha =
\sum_{i=2}^{\infty} d_{i}\alpha^{i} \in {\mathcal G}_a \subset {\mathcal R}_a, \mbox{ where } (d_{i})_{i \geq
2} \in  E(R).
\]

Since $\sum_{i=2}^{N}\varepsilon_{i}\alpha^{i}$ is a interior point
of $\mathcal{R}_a$ (item (iii) of Proposition \ref{prop12}), then
there exists a non-negative integer $M$ such that
\[
-p-q\alpha + \sum_{i=2}^{M}d_{i}\alpha^{i} = \sum_{i=2}^{\infty}
e_{i}\alpha^{i} \in \mathcal{R}_a.
\]

Since $-p-q\beta + \sum_{i=2}^{M}d_{i}\beta^{i} > 0$ then
$-p-q\beta+\sum_{i=2}^{M}d_{i}\beta^{i}= \sum_{i=l}^{K}
f_{i}\beta^{i}$ where $(f_{i})_{l \leq i \leq K} \in
 E(R)$ and $l,\ K \in
\mathbb{Z}$. Therefore,
\[
-p-q\alpha+\sum_{i=2}^{M}d_{i}\alpha^{i} =
\sum_{i=l}^{K}f_{i}\alpha^{i} = \sum_{i=2}^{\infty}e_{i}\alpha^{i}.
\]

By item (ii) of Proposition \ref{prop12},  we deduce that $e_{i}=0$
for all $i > K$ and by the same argument used in case 1, we have
that $p=q=0.$ \hspace{2.0cm} \hfill $\Box$


\begin{prop}\label{teofront2} The boundary of $\mathcal{G}_a$ satisfies the following properties:
\begin{itemize}
\item[1)] $
\partial{\mathcal{G}_{a}} = \bigcup_{u \in A} \mathcal{G}_a
\cap (\mathcal{G}_a+u) $ where $A$ is a finite set belonging to
$\mathbb{Z}+\alpha\mathbb{Z}$, whose cardinality is even and greater
than or equal to $6$ and $\{\pm 1,\pm \alpha,\pm (\alpha-1)\}
\subset A$. \label{propfront}
\item[2)] Let $z \in \partial{\mathcal{G}_{a}}$ then there exist $(\varepsilon_{i})_{i \geq 2}$ and
$(\varepsilon_{i}^{'})_{i \geq l} \in  E(G),\ l<2$ such that
$z=\sum_{i=2}^{+\infty}
\varepsilon_{i}\alpha^{i}=\sum_{i=l}^{+\infty}
\varepsilon_{i}^{'}\alpha^{i}$ and $\varepsilon_{l}^{'} \neq 0$.
\end{itemize}
\end{prop}

\dem
\begin{enumerate} \item Let $z \in
\partial{\mathcal{G}}_{a}$, then there exists a sequence $(z_{n})_{n \geq 0}$
such that
\[
\lim_{n \longrightarrow + \infty} z_{n}=z \hspace{0.3cm} \mbox{and}
\hspace{0.3cm} z_{n} \notin \mathcal{G}_a,\ \forall n \geq 0.
\]

By Theorem \ref{teo23} (a), there exists a sequence $(p_{n})_{n \geq
0}$ of elements $\mathbb{Z}+\alpha\mathbb{Z}$ such that for all $n
\geq 0,\ z_{n} \in \mathcal{G}_a+p_{n}.$ Then
$(p_{n})_{n \geq 0}$ is bounded. Since $\mathbb{Z}+\mathbb{Z}\alpha$
is a discrete group then there exists a subsequence
$(p_{k_{n}})_{n \geq 0}$  such that for all
$n,\ p_{k_{n}}=p \in \mathbb{Z}+\mathbb{Z}\alpha.$
Since $z_{k_{n}} \in \mathcal{G}_a+p$, we have  $z=\lim z_{k_{n}}
\in \mathcal{G}_a+p$.
Therefore,
\[\partial \mathcal{G}_a \subset \bigcup_{p \in \mathbb{Z}+\mathbb{Z}\alpha} \mathcal{G}_a \cap (\mathcal{G}_a+p).\]

On the other hand, if $z \in \mathcal{G}_a \cap (\mathcal{G}_a+p),\
p \in \mathbb{Z}+\mathbb{Z}\alpha \setminus \{0\}.$ Then by Theorem
\ref{teo23} (b), $z \notin int(\mathcal{G}_a)$. Therefore, $z \in
\partial \mathcal{G}_a$. Hence, $\partial \mathcal{G}_a = \bigcup_{p \in
\mathbb{Z}+ \mathbb{Z}\alpha} \mathcal{G}_a \cap (\mathcal{G}_a+p)=
\bigcup_{p \in A} \mathcal{G}_a \cap (\mathcal{G}_a+p)$ where $A=\{p
\in \mathbb{Z}+\mathbb{Z}\alpha,\ \mathcal{G}_a \cap
(\mathcal{G}_a+p) \neq \emptyset \}$.
Since $A \subset \mathbb{Z}+\mathbb{Z} \alpha \cap (\mathcal{G}_a - \mathcal{G}_a)$, we deduce that
 $A$ is finite set.  Finally, the cardinality of $A$ is
even because if $u \in A$ then $-u \in A$.

Now, we prove that $\{\pm1,\pm\alpha,\pm(\alpha-1)\} \subset A$.

In fact, it's easy to see that $-\alpha^{3}$ can be written in the
following ways:
\[
\begin{array}{ccl}
-\alpha^{3} & = & (a-1)\sum_{i=1}^{\infty}(\alpha^{4i+1}+\alpha^{4i+2})\\
&=&  \alpha+(a-2)\alpha^{3}+(a-1)\sum_{i=1}^{\infty}(\alpha^{4i}+\alpha^{4i+1}) \\
&=&  1 + (a-1) \alpha^{2} + (a-2)\alpha^{3} +
 (a-1)\sum_{i=1}^{\infty} (\alpha^{4i}+\alpha^{4i+1}).
\end{array}
\]
Hence, $-\alpha^{3} \in \mathcal{G}_a \cap (\mathcal{G}_a + \alpha)
\cap (\mathcal{G}_a+1)$. Therefore, $1$ and $\alpha$ belong to
$A$. We also show that
\[
\begin{array}{ccl}
z  & = &
\alpha-1 +(a-1) \sum_{i=1}^{\infty}(\alpha^{4i}+\alpha^{4i+1}) \\
&=&
(a-1)\sum_{i=1}^{\infty} (\alpha^{4i-2}+\alpha^{4i+1})\\
&=&  \alpha+(a-2) \alpha^{2} +
\sum_{i=1}^{\infty}(\alpha^{4i-1}+\alpha^{4i+2}).
\end{array}
\]
Then, $z \in \mathcal{G}_a \cap (\mathcal{G}_a+\alpha-1) \cap
(\mathcal{G}_a+\alpha)$. Therefore, $\alpha-1$ belongs to $A$.
\item Let $z \in \partial\mathcal{G}_a$ then \begin{equation}z=n+p\alpha+\sum_{i=2}^{+\infty} \varepsilon_{i}\alpha^{i}=
\sum_{i=2}^{\infty}
\varepsilon_{i}^{'}\alpha^{i},\label{eqteoestar}\end{equation} where
$(\varepsilon_{i})_{i \geq 2},\ (\varepsilon_{i}^{'})_{i \geq 2} \in
E(G)$  and $(n,p)  \in \mathbb{Z}^{2}\setminus \{(0,0)\}$. We have to
consider the following cases:
\begin{enumerate}
\item If  the set $\{i \geq 2,\; \varepsilon_{i} \neq 0\}$ is finite, then
$\mathcal{G}_a \cap int(\mathcal{G}_a+n+p\alpha) \neq \emptyset$,
which contradicts item b) of Theorem \ref{R}.

\item Assume that the set $\{i\geq 2,\; \varepsilon_{i} \neq 0\}$ is infinite.
Let $k \geq 2$ be an integer and
$z_{k}=n+p\alpha+\sum_{i=2}^{k}\varepsilon_{i}\alpha^{i}$.
We have that $\lim_{k \longrightarrow +
\infty} z_{k} = \sum_{i=2}^{+\infty}
\varepsilon_{i}^{'}\alpha^{i}=z$. On the other hand, there exists an integer $N \in \mathbb{N}$ such that for all $k \in
\mathbb{N}$, $\tilde{z}_{k}=n+p\beta+\sum_{i=2}^{k}\varepsilon_{i}\beta^{i} >0$. Hence  $z_{k}=\sum_{i=l_{k}}^{N_{k}}
\varepsilon_{i,k}^{''}\alpha^{i},\ (\varepsilon_{i,k}^{''} ) \in E(G)$
and $\varepsilon_{l_{k},k}^{''} > 0$. Moreover, $l_{k} < 2,$ because
otherwise $n=p=0$. On the other hand, there exists $s \in
\mathbb{Z}$ such that  $s < l_{k}$ for all integer  $ k \geq 0$
(because $(z_{k})_{k \geq 0}$ is bounded). Then, for all integer $k \geq N,\
z_{k} \in \bigcup_{t=s}^{1} \mathcal{E}_{t}$ where
$\mathcal{E}_{t}=\{\sum_{i=t}^{+\infty} \varepsilon_{i}\alpha^{i},\
(\varepsilon_{i})_{i \geq t} \in E(G),\
\varepsilon_{t}>0\}$. Since $\bigcup_{t=s}^{1} \mathcal{E}_{t}$
compact, we have $z=\lim_{k \longrightarrow +\infty} z_{k} \in
\bigcup_{t=s}^{1} \mathcal{E}_{t}$. Therefore,
$z=\sum_{i=l}^{+\infty} \varepsilon_{i}^{''}\alpha^{i},\
(\varepsilon_{i}^{''})_{l \leq i \leq +\infty} \in E(G),\
l<2$. \hfill$\Box$
\end{enumerate}
\end{enumerate}

\section{Definition of the automaton recognizing the points with at
least two expansions}
\label{def-auto}

In this section we proceed to the construction of the automaton
${\mathcal A}$ that characterize the boundary of ${\mathcal G}_a $ and ${\mathcal
R}_a $. The set of states of the automaton ${\mathcal A}$ (see Theorem
\ref{egalite}) is the set $S_{a} = \{0,\pm \alpha,\pm \alpha^2, \pm
(\alpha- \alpha^2),\;\pm (1+(a-1)\alpha^2), \pm (1+
(a-2)\alpha^2),\; \pm (1 -\alpha+ (a-1)\alpha^2),\;\pm (1-2 \alpha
+a\alpha^2)\}.$

Let $s$ and $t$ be two states. The set of edges is the set of
$(s,(e,f),t)\in S \times \{ 0,1,\ldots, a-1 \}^2 \times S$ satisfying
$t = \frac {s}{\alpha} + (e-f) \alpha^{2}$. The set of initial
states is $\{ (0,(0,0),0) \}$.

Let us explain the behaviour of this automaton. Let $\varepsilon =
(\varepsilon_{i})_{i\geq l}$ and $\varepsilon' =
(\varepsilon'_{i})_{i\geq l}$ belonging $E(R)$ (resp. E(G)),
$x=\sum_{i=l}^{\infty}\varepsilon_{i}\alpha^{i}$ and
$y=\sum_{i=l}^{\infty}\varepsilon'_{i}\alpha^{i}$. Suppose $x=y$.
For all $k\geq l$ we put

\begin{align}
\label{def-Ak}
A_k (\varepsilon ,\varepsilon')
=
\alpha^{-k+2}\sum_{i=l}^{k}(\varepsilon_{i}-\varepsilon'_{i})\alpha^{i}.
\end{align}

In the following we prove that all the $A_k$, $k\in \mathbb{N}$,
belong to $S$. Clearly, for all $k\geq l$,

\begin{equation}\label{relaq}
A_{k+1} (\varepsilon ,\varepsilon') = \frac {A_{k}(\varepsilon
,\varepsilon')}{\alpha} + (\varepsilon_{k+1}- \varepsilon'_{k+1})
\alpha^{2} .
\end{equation}

Let $s$ be the smallest integer such that $\varepsilon_{s} \ne
\varepsilon'_{s}$. Hence $A_{i} (\varepsilon ,\varepsilon') = 0$ for
$i \in \{l, \ldots, s-1\}$. Suppose  $\varepsilon_{s}
>\varepsilon'_{s}$. Then, $A_{s}= (\varepsilon'_{s}
-\varepsilon_{s})\alpha^{2}= \alpha^{2}$. From (\ref{relaq}) we
deduce $ A_{s+1}(\varepsilon ,\varepsilon') = \alpha +
(\varepsilon_{s+1}- \varepsilon'_{s+1}) \alpha^{2}$ which should
belong to $S_{a,b}$. Hence $A_{s+1} (\varepsilon ,\varepsilon')=
\alpha $ if $\varepsilon _{s+1}= \varepsilon'_{s+1}$ or  $A_{s+1}
(\varepsilon ,\varepsilon')= \alpha- \alpha^{2} $ if $(\varepsilon
_{s+1}, \varepsilon'_{s+1})= (t_1, t_1 +1)$, where $ 0 \leq t_1 \leq
a-2$. Continuing by the same way and using the fact that
 the set of states   $S$ is finite, we obtain a
finite state automaton.
\begin{obs}
The idea of using finite state automaton to recognize points that
have at least 2 $\alpha$- expansions is old. It was done in the case
of $\alpha = 1/ \gamma$ where $\gamma >1$ is a Pisot number and the
digits belong to a finite set of integer numbers  by Frougny in
\cite{frougny1}. In \cite{Th} Thurston proved the same result in the case where
$\beta$ is a Pisot complex numbers and the digits are in a finite
subset of algebraic integers in $\mathbb{Q} (\gamma) $ (see also \cite{HZ98}, \cite{ali05}).
 The difficulty remains in the fact that it is
not easy to find exactly the set of states. The classical method
uses the modulus of $\alpha$. In this work, we give a method which
does not use the modulus of $\alpha$, with this we could find all
the states for the automata associated to a class of cubic Pisot
unit numbers.

\end{obs}

\subsection{Characterization of the points with two expansions}
\label{sec-characterization}

Let $\varepsilon = (\varepsilon_{i})_{i\geq l}$ and $\varepsilon' =
(\varepsilon'_{i})_{i\geq l}$ in $E(R)$ (resp. $E(G)$) where $l \in
\mathbb{Z},\; x=\sum_{i=l}^{\infty}\varepsilon_{i}\alpha^{i}$ and
$y=\sum_{i=l}^{\infty}\varepsilon'_{i}\alpha^{i}$. Suppose that
$x=y$. For all $k\geq l$ put
\begin{align}
\label{def-Ak} A_k (\varepsilon ,\varepsilon')=A_{k} =
\alpha^{-k+2}\sum_{i=l}^{k}(\varepsilon_{i}-\varepsilon'_{i})\alpha^{i}.
\end{align}
Then,
\begin{eqnarray}
\label{rel-Ak} A_{k+1} = \frac {A_{k}}{\alpha} + (\varepsilon_{k+1}-
\varepsilon'_{k+1}) \alpha^{2}.
\end{eqnarray}
For all $\varepsilon = (\varepsilon_{i})_{i \geq l}$ and
$\varepsilon'=(\varepsilon'_{i})_{i \geq l}$ in $E(R)$ (resp. E(G)),
let
$$ S(\varepsilon, \varepsilon')=\{ A_k (\varepsilon,\varepsilon')
; k\geq l \} = \left\{
\alpha^{-k+2}\sum_{i=l}^{k}(\varepsilon_{i}-\varepsilon'_{i})\alpha^{i}
; k \geq l \right\} .
$$

\begin{teo}
\label{egalite} Let
$x=\sum_{i=l}^{\infty}\varepsilon_{i}\alpha^{i}$,
$y=\sum_{i=l}^{\infty}\varepsilon'_{i}\alpha^{i}$, where
 $\varepsilon = (\varepsilon_{i})_{i \geq l}$ and
$\varepsilon' =(\varepsilon'_{i})_{i \geq l}$ in $E(R)$ (resp. $E(G)$). Thus, $x=y$ if only if $S(\varepsilon, \varepsilon')$
is finite. Moreover, $ S(\varepsilon ,\varepsilon') \subset S_{a} =
\{0,\pm \alpha,\pm \alpha^2, \pm (\alpha- \alpha^2),\;\pm
(1+(a-1)\alpha^2), \pm (1+ (a-2)\alpha^2),\; \pm (1 -\alpha+
(a-1)\alpha^2),\;\pm (1-2 \alpha +a\alpha^2)\}.$ And
$$
S_{a}=\bigcup_{(\varepsilon ,\varepsilon')\in \Delta  } S
(\varepsilon ,\varepsilon').
$$
where $\Delta = \left\{ \left((\varepsilon_i)_{i\geq l},
(\varepsilon'_i)_{i\geq l}\right) \in X \times X ;
\sum_{i=l}^{\infty}\varepsilon_{i}\alpha^{i}=\sum_{i=l}^{\infty}\varepsilon'_{i}\alpha^{i}
\right\}$,
where $X= E(R)$  (resp. $X= E(G)$).
\end{teo}

\thinlines    
\dem If $S(\varepsilon,\varepsilon')$ is finite, then since $0 < |\alpha| < 1,$ we have
$\sum_{i=l}^{+\infty}\varepsilon_{i}\alpha^{i}=\sum_{i=l}^{+\infty}\varepsilon_{i}'\alpha^{i}$.

Let $x=\sum_{i=l}^{\infty}\varepsilon_{i}\alpha^{i}$ and
$y=\sum_{i=l}^{\infty}\varepsilon_{i}'\alpha^{i}$ with $\varepsilon
= (\varepsilon_{i})_{i \geq l}$ and
$\varepsilon'=(\varepsilon'_{i})_{i \geq l}$ in $E(R)$  (resp. $E(G)$).
Assume $x=y$. Then for all $k\geq l$, we have
\begin{eqnarray}
 \label{Ak}
A_k =
 \sum_{i=k+1}^{\infty}(\varepsilon'_{i}-\varepsilon_{i})\alpha^{i-k+2}
 =
 \sum_{i=3}^{\infty}(\varepsilon'_{i+k-2}-\varepsilon_{i+k-2})\alpha^{i} \in  \alpha {\mathcal R}- \alpha{\mathcal
 R}\
 (\mbox {resp. } \alpha {\mathcal G}- \alpha {\mathcal G}).
\end{eqnarray}
Let $k\geq l$ and assume that $A_k \not = 0$. By \eqref{def-Ak} and
the fact that $\alpha$ is an algebraic integer of degree $3$, we
deduce that
\begin{align}
\label{Ak-alg} A_k= n_{k} \alpha^{2}+ p_{k} \alpha+ q_{k}, \mbox{
where } n_{k},\ p_{k},\ q_{k} \in \mathbb{Z}.
\end{align}
Put $\widetilde{A_{k}}=n_{k}\beta^{2}+p_{k}\beta + q_{k}$. Since
$\widetilde{A_{k}}$ or $-\widetilde{A_{k}}$ belong to
$\mathbb{Z}[\beta] \cap \mathbb{R}^{+}$, we deduce according to item
(iv) of Proposition \ref{prop12}, that there exists a sequence
$(c_{i})_{s_k \leq i \leq m_k} \in E(R)$   such that
$c_{m_{k}}>0$ and
\begin{align}
 \label{Ak-alg2}
\widetilde{A_{k}}=n_{k} \beta^{2}+ p_{k} \beta+ q_{k}  = \pm
\sum_{i=s_{k}}^{m_{k}}c_{i}\beta^{i}.
\end{align}
Assume that
$\widetilde{A_{k}}=\sum_{i=s_{k}}^{m_{k}}c_{i}\beta^{i}$. By using  \eqref{def-Ak},
\eqref{Ak-alg}, \eqref{Ak-alg2} and the fact that
$\beta $ and $\alpha$ are algebraic conjugates, we get
\begin{align}
\label{eq-ali} \beta^{-k+2}\sum_{i=l}^{k} \varepsilon_{i} \beta^{i}
& = \beta^{-k+2}\sum_{i=l}^{k} \varepsilon'_{i}\beta^{i} +
\sum_{i=s_{k}}^{m_{k}}c_{i}\beta^{i}.
\end{align}
According to (vi) of Proposition \ref{prop12}, we deduce that
$\beta^{-k+2}\sum_{i=l}^{k}\varepsilon_{i} \beta^{i} < \beta^3$.
Consequently $m_k \leq 2$. Putting $c_i = 0$ for all $i>m_k$, we
have $n_{k} \beta^{2}+ p_{k} \beta+ q_{k}  =
\sum_{i=s_{k}}^{2}c_{i}\beta^{i}$. Since $\beta$ is a Pisot number,
using Proposition [2] in \cite{SF}, we deduce that there
exists an integer $s=s(a)$ such that $s \leq s_{k}$. Therefore,
\begin{align}
\label{Akc} A_k  =  \sum_{i=s}^{2} c_{i}\alpha^{i}.
\end{align}
Then \begin{eqnarray}S_{a} \subset \{\sum_{i=s}^{2} c_{i}\alpha^{i},
(c_{i})_{s \leq i \leq 2} \in E(R)  \}.
\label{eq50}\end{eqnarray} Note that if
\begin{equation}
A_{k}=\sum_{i=s}^{2}c_{i}\alpha^{i} \mbox{ then }
\widetilde{A}_{k}=\sum_{i=s}^{2}c_{i}\beta^{i}< \beta^{3}.
\label{eqbeta3}
\end{equation}
To prove that $S_{a}=  \{0,\pm \alpha,\pm \alpha^2, \pm (\alpha-
\alpha^2),\;\pm (1+(a-1)\alpha^2), \pm (1+ (a-2)\alpha^2),\; \pm (1
-\alpha+ (a-1)\alpha^2),\;\pm (1-2 \alpha +a\alpha^2)\},$ we need
the following important result:
\begin{prop}
Let $n,p,q \in \mathbb{Z}$ and $z=n+p\alpha+q\alpha^{2}$. If $z\in
S_{a}$ then $|n| \leq 1$. \label{lem3cap}
\end{prop}
\dem Let $S_{a}=\{A_{k}=n_{k}+p_{k}\alpha+q_{k}\alpha^{2},\ k \geq
0\}$ and $n=max \{|n_{k}|,\ k \geq 0\}$. Suppose $n \geq 2$ and let
$k \in \mathbb{N}$ such that $A_{k}=n+p\alpha+q\alpha^{2},\ p,q \in
\mathbb{Z}$. Thus by (\ref{rel-Ak}),
$A_{k+1}=\frac{n}{\alpha}+p+q\alpha+ d\alpha^{2},|d| \leq a-1$, and
using that $\alpha^{-1}=\alpha^2-a\alpha+1$ we have
\[A_{k+1}=(p+n)+(q-na)\alpha+(d+n)\alpha^{2}.\]
Therefore $|p+n|\leq n$ and then $-2n \leq p \leq 0$.\\ On the other
hand, $ A_{k+2}=(q-na+p+n)+(d+n-a(p+n))\alpha+(f+p+n)\alpha^{2},\
\mbox{ where } |f| \leq a-1$, and by the same reason $n(a-2)-p \leq
q$. Then we have
\begin{equation}
\label{eqeq1001}
\begin{array}{l}\widetilde{A}_{k}=n+p\beta+q\beta^{2}=1-\beta+a \beta^{2}
+(n-1)+ (p+1)\beta+(q-a)\beta^{2}\\ \geq \beta^{3}+(n-1)+(p+1)\beta
+ (n(a-2)-p-a)\beta^{2},\end{array}
\end{equation}
and since $n \geq 2$ then $n(a-2)-p-a \geq a-4-p$. Thus we conclude
by (\ref{eqeq1001}) that
\begin{equation} \widetilde{A}_{k} \geq
\beta^{3}+(n-1)+(p+1)\beta+ (a-4-p)\beta^{2}. \label{eqeq02}
\end{equation}
We have to analyze all the following cases.
\begin{enumerate}
\item[1-] $a\geq4$ or ($a=3$ and $-2n<p<0$).\\ Here we have $a-4-p > 0$ and
 follows by (\ref{eqeq02}) that $$\widetilde{A}_{k} \geq
\beta^{3}+(n-1)+(p+1)\beta+ (a-4-p)\beta^{2}\geq
\beta^{3}+(n-1)+(p+1+a-4-p)\beta\geq\beta^3.$$ This cannot occur
because of (\ref{eqbeta3}).


\item[2-] If $a=3$ and $p=0$.\\ Here we have $A_{k+3}=(d+q-4n)+(f+7n-3q)\alpha+(g+q-2n)\alpha^2,\ |g| \leq a-1=2,$
and since $|d+q-4n|\leq n$, $|d|\leq 2$ then $-n\leq
d+q-4n\leq2+q-4n$. Thus
\begin{equation}
1 \leq 3n-5\leq q-3. \label{eqeq03}
\end{equation}
By (\ref{eqeq1001}) and $n\geq 2$, we have
\[
\widetilde{A}_{k} \geq \beta^{3} +(n-1) + \beta + (q-3)\beta^{2} >
\beta^{3}.
\]

\item[3-] If $a=2$ and $p \leq -3$.\\
Again $a-4-p>0$, and by (\ref{eqeq02}) $\widetilde{A}_{k} \geq
\beta^{3}+(n-1)+(p+1)\beta+(-2-p)\beta^2$. We have
$\beta^2=2\beta-1+\beta^{-1}$, $\beta>1$ and therefore
\[
\begin{array}{lll}
\widetilde{A}_{k} &\geq& \beta^{3}+(-2-p)\beta^{-1} + (n+1+p)+(-3-p)\beta \\

 &>& \beta^{3} + (-2-p)\beta^{-1}+(n-2) \geq \beta^{3}.
\end{array}
\]
Therefore if $A_{k}=n+p\alpha+q\alpha^2$ then
\begin{equation}\label{estrela}
a=2 \mbox{ and } p\in\{-2,-1,0\}.
\end{equation}
\item[4-] If $a=2$ and $-2 \leq p \leq 0$.\\ Here the possibilities are $A_{k}=n+q\alpha^{2}$ or $A_{k}=n-\alpha+q\alpha^{2}$ or
$A_{k}=n-2\alpha+q\alpha^{2}$.\\
$\bullet$ If $A_{k}=n+q\alpha^{2}.$
\\We have $A_{k+1}=n+(q-2n)\alpha+(d+n)\alpha^{2}$
and, using the previous cases, we get
$q=2n$ or $q=2n-1$ or $q=2n-2$ and thus $\widetilde{A}_{k} = n
+ q \beta^{2}> \beta^3.$
\\ $\bullet$ $A_{k}=n-\alpha+q\alpha^{2}.$\\
Then
\[
\begin{array}{l}
A_{k+1} = (n-1) + (q-2n) \alpha + (d+n) \alpha^{2} \mbox{ and } \\
A_{k+2} = (q-n-1) + (d+n-2(n-1))\alpha + (f+n-1)\alpha^{2}.
\end{array}
\]
Since $-n \leq q-n-1 \leq n$ then $ 1 \leq q$.\\ If $q \geq 2=a$
then
\[
\widetilde{A}_{k}=n-\beta+q\beta^{2} \geq n -\beta + a \beta^{2}
\geq \beta^{3}.
\]
If $q=1$, then $A_{k+2}=-n+(d+n-2(n-1))\alpha+(f+n-1)\alpha^{2}$.
Since $A_{k+2} \in S_{a}$, then
 $-A_{k+2} \in S_{a}$. Therefore
\[
-A_{k+2}=n+(2(n-1)-n-d)\alpha+(-f-n+1)\alpha^{2}.
\]

Using (\ref{estrela}) we get $2(n-1)-n-d  \in \{-1,-2\}$.

If $2(n-1)-n-d  =-2$, then $d = n \leq a-1=1$. That is absurd, since
$n \geq 2$.\\ If $2(n-1)-n-d  =-1$, then $d = n-1 \leq 1$, thus $n
\leq 2$.

Hence $n=2$ and thus $A_{k}=2-\alpha+\alpha^{2}.$\\ If
$A_{k}=2-\alpha+\alpha^{2}$, then
\[
\begin{array}{lll}
A_{k+1} &=& 1 +(1-2a)\alpha+(d+2)\alpha^{2} \\
&=& 1 -3\alpha+(d+2)\alpha^{2},
\end{array}
\]
and
\[
\begin{array}{lll}
A_{k+2} &=& -2 +(d+2-a)\alpha+(f+1)\alpha^{2} \\
&=& -2 +d\alpha+(f+1)\alpha^{2}.\
\end{array}
\]
We know that $A_{k+2}=-2+\alpha-\alpha^{2}$ or
$A_{k+2}=-2+2\alpha-\alpha^{2}$. In the first case, we obtain $d=1$
and $f=-2$, which is impossible because
 $|f| \leq a-1=1$. In the second case, we get  $d=2$ which is also absurd.\\
%

$\bullet$ $A_{k}=n-2 \alpha+q \alpha^{2}$.\\ Then $A_{k+1}=
(n-2)+(q-2n)\alpha+(d+n)\alpha^{2}$ and
\[ A_{k+2} = (q-2n+n-2)+(d+n-2(n-2))\alpha + (f+n-2)\alpha^{2}, |f| \leq 1.\]

Since $-n \leq q-2n+n-2 \leq n$ then  $2
\leq q$.\\
If $q \geq 3$, then
\[
\begin{array}{lll}
\widetilde{A}_{k} &=& n-2\beta+q\beta^{2} \geq n-2\beta+(a+1)\beta^{2} \\
&=& \beta^{3} + (\beta^{2}-\beta) > \beta^{3}.

\end{array}
\]
If $q=2$, then $A_{k}=n-2\alpha+2\alpha^{2}$ and $A_{k+2}= - n +
(d-n+4)\alpha + (f+n-2)\alpha^{2}$. Since $-A_{k+2}\in S$ we have
$d-n+4=2$ and therefore
$d=n-2 \leq 1$ and $n \leq 3$.\\
If $n=2$, then $A_{k}= 2 - 2\alpha + 2\alpha^{2}$, $A_{k+1}=
(2-2a)\alpha+(d+2)\alpha^{2}=-2\alpha + (d+2)\alpha^{2}.$ Then $
A_{k+2}=-2+(d+2)\alpha+f \alpha^{2}$,  hence $ d=0$   and $f=2$.
That is impossible.\\
If $n=3$, then $A_{k}= 3 - 2\alpha + 2\alpha^{2}$, $A_{k+1}=
1+(2-3a)\alpha+(d+3)\alpha^{2}=1-4\alpha + (d+3)\alpha^{2}$, then
$A_{k+2}=-3+(d+3-a)\alpha+(f+1)\alpha^{2}$. Hence $d+3-a=2$, thus
$d=a-1=1$ and $f+1=-2$. Therefore
   $f=-3$. That is an absurd.
\end{enumerate}
\hfill$\Box$


Now we will prove that\\ $S_{a}=  \{0,\pm \alpha,\pm \alpha^2, \pm
(\alpha- \alpha^2),\;\pm (1+(a-1)\alpha^2), \pm (1+
(a-2)\alpha^2),\; \pm (1 -\alpha+ (a-1)\alpha^2),\;\pm (1-2 \alpha
+a\alpha^2)\}.$\\ If $A_{k}=n+p\alpha+q\alpha^{2} \in S_{a},$ then
by proposition \ref{lem3cap} we have $|n| \leq 1$. Consider the
following cases:
\begin{enumerate}
\item[1-] $n=1$, then we have
$A_{k}=1+p\alpha+q\alpha^{2}$ and
$A_{k+1}=(p+1)+(q-a)\alpha+(d+1)\alpha^{2}$, where $|d| \leq a-1$.
Hence by Proposition \ref{lem3cap}, $p \in \{-2,-1,0\}$.

$\bullet$ If $p=0$ then $A_{k+1}=1+(q-a)\alpha+(d+1)\alpha^{2}$.
Thus we have $a-2 \leq q\leq a$. If $q=a$ then
$\widetilde{A}_{k}=1+a \beta^{2}> \beta^{3}$, which is impossible
because of (\ref{eqbeta3}). Hence we have the states
$A_{k}=1+(a-1)\alpha^{2}$ or $A_{k}=1+(a-2)\alpha^{2}$. \\
$\bullet$ If $p=-1$ then $A_{k+1}=(q-a)\alpha + (d+1)\alpha^{2}$ and
$A_{k+2}=(q-a)+ (d+1)\alpha + e \alpha^{2},\ |e| \leq a-1.$ Hence we
have $-1 \leq q-a \leq 1$. For the cases $q=a$ or $q=a+1$ we have
$\widetilde{A}_{k} \geq \beta^{3}.$ Hence we get the
state $A_{k}=1-\alpha+(a-1)\alpha^{2}$.\\
$\bullet$ If $p=-2$, then $A_{k+1}=-1+(q-a)\alpha + (d+1)\alpha^{2}$
and $A_{k+2}=(q-a-1)+ (d+1+a)\alpha + e \alpha^{2},\ |e| \leq a-1.$
Hence $-1 \leq q-a-1 \leq 1$. If $q=a$ we get the state
$A_{k}=1-2\alpha+a\alpha^{2}$. If $q=a+1$ or $q=a+2$, then
$\widetilde{A}_{k}> \beta^{3}$.

\item[2-] $n=0$, then we have $A_{k}=p\alpha+q\alpha^{2}$ and
$A_{k+1}=p+q\alpha+d\alpha^{2}$. Hence  by Proposition
\ref{lem3cap},
$p\in\{-1,0,1\}$.\\
$\bullet$ If  $p=0$ then $A_{k+1}=p+q\alpha+d\alpha^{2}$,
$A_{k+2}=q+d\alpha+e\alpha^{2}$ and by Proposition \ref{lem3cap}, we
deduce that $q=0, \pm 1$. Therefore, $A_{k}=0,\ \pm \alpha^{2}$.\\
$\bullet$ If $p=1$ then $A_{k+1}=1+q\alpha + d\alpha^{2}$ and
$q\in\{0,-1,-2\}$. Therefore,
 $A_{k}=\alpha,\ \alpha-\alpha^{2},\
\alpha-2\alpha^{2}$. If $A_{k}=\alpha-2\alpha^{2}$ is a state then
$A_{k+1}=1-2\alpha+d\alpha^{2} $is also a state. But as noted
earlier $A_{k+1}=1-2\alpha+d\alpha^{2}=1-2\alpha+a\alpha^{2}$ and,
therefore $d=a$, which is impossible since  $d \leq a-1$.\\
$\bullet$ If $p=-1$. Using the same ideas of previous case, we
obtain the states $A_{k}=-\alpha,\ -\alpha+\alpha^{2}$. \hfill$\Box$


\begin{picture}(518,505)
\thinlines    
              \put(264,356){\scriptsize (a-1,0)}
              \put(91,360){\scriptsize (0,a-1)}
              \put(82,410){\scriptsize ($\varepsilon+1,\varepsilon$)}
              \put(245,250){\scriptsize ($\varepsilon+a-2,\varepsilon$)}
              \put(118,250){\scriptsize ($\varepsilon,\varepsilon+a-2$)}
              \put(72,313){\scriptsize ($\varepsilon,\varepsilon+a-2$)}
              \put(282,313){\scriptsize ($\varepsilon+a-2,\varepsilon$)}
               \put(288,413){\scriptsize ($\varepsilon$,$\varepsilon+1$)}
              \put(129,287){\scriptsize (0,a-1)}
              \put(207,176){\scriptsize (a-1,0)}
              \put(240,168){\scriptsize (a-1,0)}
              \put(172,168){\scriptsize (0,a-1)}
              \put(193,198){\scriptsize (0,a-1)}
              \put(252,286){\scriptsize (a-1,0)}
              \put(269,472){\scriptsize (1,0)}
              \put(125,472){\scriptsize (0,1)}
              \put(155,380){\scriptsize ($\varepsilon$,$\varepsilon$)}
              \put(224,380){\scriptsize ($\varepsilon$,$\varepsilon$)}
              \put(165,340){\scriptsize ($\varepsilon$,$\varepsilon$)}
              \put(157,315){\scriptsize ($\varepsilon+1$,$\varepsilon$)}
              \put(217,315){\scriptsize ($\varepsilon$,$\varepsilon+1$)}
              \put(270,333){\scriptsize $1-\alpha+(a-1)\alpha^{2}$}
              \put(90,333){\scriptsize $-1+\alpha-(a-1)\alpha$}
              \put(280,293){\scriptsize $1+(a-1)\alpha^{2}$}
              \put(70,293){\scriptsize $-1-(a-1)\alpha^{2}$}
              \put(235,293){\scriptsize $\alpha$}
              \put(159,293){\scriptsize $-\alpha$}
              \put(235,203){\scriptsize $1+(a-2)\alpha^{2}$}
              \put(223,143){\scriptsize $1-2\alpha+a\alpha^{2}$}
              \put(128,143){\scriptsize $-1+2\alpha-a\alpha^{2}$}
              \put(128,203){\scriptsize $-1-(a-2)\alpha^{2}$}
             \put(268,383){\scriptsize $\alpha-\alpha^{2}$}
              \put(100,385){\scriptsize $-\alpha+\alpha^{2}$}
              \put(109,438){\scriptsize $-\alpha^2$}
              \put(279,438){\scriptsize $\alpha^2$}
              \put(200,489){\scriptsize 0}
              \put(288,378){\vector(0,-1){30}}
              \put(289,305){\vector(0,1){20}}
              \put(115,378){\vector(0,-1){30}}
              \put(115,430){\vector(0,-1){30}}
              \put(279,430){\vector(0,-1){30}}
              \put(251,295){\vector(1,0){20}}
              \put(234,445){\vector(0,-1){140}}
              \put(164,435){\vector(0,-1){130}}
             \put(234,445){\line(1,0){30}}
              \put(164,435){\line(-1,0){30}}
              \put(151,295){\vector(-1,0){22}}
               \put(165,284){\vector(0,-1){70}}
              \put(240,284){\vector(0,-1){70}}
              \put(120,305){\vector(0,1){20}}
              \put(238,195){\vector(0,-1){40}}
              \put(220,155){\vector(-2,3){36}}
              \put(188,152){\vector(2,3){36}}
              \put(267,325){\vector(-2,1){120}}
               \put(156,330){\vector(2,1){100}}
              \put(170,195){\vector(0,-1){40}}
              \put(229,490){\vector(2,-1){65}}
              \put(262,332){\vector(-2,-1){78}}
              \put(155,325){\vector(2,-1){68}}
              \put(180,490){\vector(-2,-1){65}}
              \put(269,325){\framebox(68,18){}}
              \put(90,325){\framebox(66,18){}}
              \put(90,379){\framebox(47,18){}}
              \put(274,285){\framebox(60,18){}}
              \put(226,285){\framebox(25,18){}}
              \put(152,285){\framebox(25,18){}}
              \put(68,285){\framebox(60,18){}}
              \put(100,432){\framebox(30,18){}}
              \put(265,432){\framebox(30,18){}}
              \put(228,196){\framebox(60,18){}}
              \put(220,135){\framebox(60,18){}}
              \put(124,135){\framebox(60,18){}}
              \put(124,196){\framebox(60,18){}}
              \put(260,377){\framebox(47,18){}}
              \put(180,484){\framebox(47,18){}}
\vspace{-5cm}\hspace{5cm} Automaton ${\mathcal A}$
\end{picture}
\ \\

As a Corollary we obtain the following result:
\begin{teo}\label{teo6regioes} For all $a \geq 2$ we have \[
\partial{\mathcal{G}}_{a} = \bigcup_{u \in B} \mathcal{G}_{a} \cap
(\mathcal{G}_{a}+u)
\] where $B = \{\pm 1,\pm \alpha, \pm (\alpha-1)\}$.
\end{teo}
\demT  According to item (1) of Proposition \ref{teofront2}, we have
that $\bigcup_{u \in B} \mathcal{G}_a \cap (\mathcal{G}_a+u) \subset
\partial \mathcal{G}_a$.
If $z \in \partial{\mathcal{G}}_{a}$ then by item (2) of  Proposition
\ref{teofront2} we have
\begin{equation}
z=\sum_{i=2}^{\infty} \varepsilon_{i} \alpha^{i} =
\sum_{i=l}^{\infty} \varepsilon_{i}^{'}\alpha^{i} \mbox{ where } l<2
\mbox{ and } \varepsilon_{l}^{'} \neq 0.  \label{eqteo01}
\end{equation}

Then $r=(\varepsilon_{l}^{'},0)
(\varepsilon_{l+1}^{'},0) \ldots
(\varepsilon_{1}^{'},0)(\varepsilon_{2}^{'}, \varepsilon_{2}) \ldots
,\; l <2$ is a path in the automaton ${\mathcal A}$ starting from the initial state. Therefore,
\begin{equation}
r=(1,0)(x_0,x_0+1)(a-1,0) \ldots \label{eeqeq2002}
\end{equation}
or
\begin{equation}
r = (1,0)(x_0,x_0)(a-1,0)(x_{1}+a-2,x_{1}) (x_{2},x_{2})(0,a-1) \ldots
\label{eeqeq2003}
\end{equation}
or
\begin{equation}
r = (1,0)(x_0,x_0)(a-1,0)(x_1+a-2,x_1) (x_2,x_2+1) \ldots
\label{eeqeq2004}
\end{equation}
or
\begin{equation}
r = (1,0)(x_0,x_0)(x_1+a-2,x_1)yyy \mbox{ where } y=(a-1,0)
(0,a-1)(0,a-1)(a-1,0)  \label{eeqeq2005}
\end{equation}
where $x_0,x_1,x_2 \in \{0,1,\ldots, a-1\}$.
\begin{enumerate}
\item[1-] If $r=(1,0)(x_0,x_0+1)(a-1,0) \ldots$ then
by (\ref{eqteo01}) we have $l=1$ and
\begin{equation}
z=(x_0+1)\alpha^{2}+\sum_{i=4}^{\infty} \varepsilon_{i}\alpha^{i} =
\alpha +
x_0\alpha^{2}+(a-1)\alpha^{3}+\sum_{i=4}^{\infty}\varepsilon_{i}^{'}\alpha^{i}.
\label{eeqeq25}
\end{equation}
Therefore $z \in \mathcal{G}_a \cap (\mathcal{G}_a+\alpha).$
\item[2-]If $r = (1,0)(x_0,x_0)(a-1,0)(x_{1}+a-2,x_{1}) (x_{2},x_{2})(0,a-1) \ldots$, then  $x_1=0$ and we have  the following possibilities:
\begin{enumerate}
\item[2.1-] If $x_0>0$ then $l=1$ and
 $z \in \mathcal{G}_a \cap (\mathcal{G}_a + \alpha)$.
\
\item[2.2-] If $x_0=0$ and $x_2 > 0$ then $l \in \{-2, -1,0,1\}$.
If $l=0$ or $1$, then $z \in \mathcal{G}_a \cap (\mathcal{G}_a+\alpha^l)$.

If $l=-1$, then  we have
by equation (\ref{eqteo01})
$$
\begin{array}{ccl}
z&=& x_2\alpha^{3} + (a-1)\alpha^{4}+ \sum_{i=4}^{\infty}
\varepsilon_{i}\alpha^{i}\\
 &=& \frac{1}{\alpha} +(a-1) \alpha+
(a-2)\alpha^2+\sum_{i=3}^{\infty}\varepsilon_{i}^{'}\alpha^{i}\\
&=&  1-\alpha +
(a-1)\alpha^{2}+\sum_{i=3}^{\infty}\varepsilon_{i}^{'}\alpha^{i}.
\end{array}
$$
Hence $z \in \mathcal{G}_a \cap (\mathcal{G}_a-\alpha+1)$.

If $l=-2$, then  we have
by equation (\ref{eqteo01})
\[
\begin{array}{ccl}
z&=& x_2\alpha^{2} + (a-1)\alpha^{3}+ \sum_{i=4}^{\infty}
\varepsilon_{i}\alpha^{i}
\\ &=& \frac{1}{\alpha^{2}} +(a-1)+
(a-2)\alpha+x_2\alpha^{2}+\sum_{i=3}^{\infty}\varepsilon_{i}^{'}\alpha^{i}\\
&=&  -\alpha +
(x_2+1)\alpha^{2}+\sum_{i=3}^{\infty}\varepsilon_{i}^{'}\alpha^{i}.
\end{array}
\]
Therefore, $z \in \mathcal{G}_a \cap (\mathcal{G}_a-\alpha)$.

\item[2.3-] If $x_0=x_2=0$ then $l \in \{-3,-2,-1,0,1\}$.
If $l=-3$, then
\[
z=(a-1)\alpha^{2}+\sum_{i=3}^{\infty}\varepsilon_{i}\alpha^{i}=
\frac{1}{\alpha^{3}}+\frac{(a-1)}{\alpha} +(a-2)+\sum_{i=2}^{\infty}
\varepsilon_{i}^{'}\alpha^{i}.
\]
Since $\frac{1}{\alpha^{3}}=(1-a)+\alpha+\frac{(1-a)}{\alpha}$, we
have $z \in \mathcal{G}_a \cap (\mathcal{G}_a+(\alpha-1))$.
\end{enumerate}
\item[3-] If $r=(1,0)(x,x)(a-1,0)(x_{1}+a-2,x_{1})
(x_{2},x_{2}+1) \ldots$ , then we can prove by the same way than the previous cases that
 $z \in \mathcal{G}_a\cap(\mathcal{G}_a+\alpha)$ or $z \in
\mathcal{G}_a \cap (\mathcal{G}_a+(1-\alpha))$ or $z \in
\mathcal{G}_a\cap (\mathcal{G}_a-\alpha).$
\item[4] Finally considering the path: $r=(1,0)(x,x)(x_1+a-2,x_1)
sss$ where \\ $s=(a-1,0)(0,a-1)(0,a-1)(a-1,0)$ and analyzing all possible
cases, we have: $z \in \mathcal{G}_a\cap(\mathcal{G}_a+\alpha)$ or $z
\in \mathcal{G}_a\cap(\mathcal{G}_a+1)$ or $z \in
\mathcal{G}_a\cap(\mathcal{G}_a-1)$.
\hfill$\Box$
\end{enumerate}

\section { Rauzy fractal ${\mathcal R}_a$}

Now, let us consider the classical Rauzy fractal ${\mathcal R}_a$
associated to the Pisot unit number $\beta >1$ satisfying $\beta^3
-a \beta^2+\beta-1$. We have
$$
{\mathcal R}_{a}=
\{\sum_{i=2}^{\infty}
\varepsilon_{i}\alpha^{i} \mid \forall i \geq 2,\
\varepsilon_{i}=0,1,\ldots,a-1,\ \varepsilon_{i}
\varepsilon_{i-1}\varepsilon_{i-2}\varepsilon_{i-3} <_{lex}
(a-1)(a-1)01, \mbox{ where } \varepsilon_{1}=\varepsilon_{0}= \varepsilon_{-1}=0 \}.
$$

As we mentioned before, the set
${\mathcal R}_{a}$ is a compact, connected subset of $\mathbb{C}$ with interior simply connected. Moreover, ${\mathcal R}_{a}$ induces a periodic tiling of the plane $\mathbb{C}$.


\begin{prop}\label{prop35}
${\mathcal R}_{a}$ induces a periodic tiling of the plane $\mathbb{C}$ modulo
$\mathbb{Z}u+\mathbb{Z} \alpha u$ where $u= \alpha-1$. Moreover
 $\partial{\mathcal R}_a = \bigcup_{v \in B} \mathcal{R}_a \cap (\mathcal{R}_a + v)$, where
  $B= \{\pm u, \pm  \alpha u, (1+\alpha) u, (1-\alpha) u \}$
and $ \mathcal{R}_a \cap (\mathcal{R}_a +  (1+ \alpha) u))= \{-1\},\;
 \mathcal{R}_a \cap (\mathcal{R}_a +  ( \alpha-1) u)= \{-\alpha\}$.

\end{prop}

{\bf Proof}: Consider  the sequence $ (R'_{n})_{n \geq 0}$ by
$R'_{0}=0,\ R'_{1}=0,\ R'_{2}=1,\ R'_{3}=a,\ R'_{4}=a^2,\;
R'_{n+3}=a R'_{n+2}-R'_{n+1}+R'_{n},\ n \geq 2$. Then $R_{n}=
R'_{n+2}$ for all integer $n \in \mathbb{N}$.

On the other hand, we can prove by induction on $n$ that $G'_n=
R'_{n} - G'_{n-2}$ for all integer  $n \geq 2$. Thus, using item (vii) of Proposition \ref{prop12}, we have
$$\alpha^n = R'_n \alpha^2 - G'_{n-2} (\alpha^2- \alpha) - G'_{n-1}(\alpha-1),\; \forall n \geq 2.$$
Thus $(\varepsilon_i)_{2 \leq i \leq N} \in E(R)$, we have
$$\sum_{i=2}^{N}\varepsilon_i \alpha^i= n \alpha^2 + p_n (\alpha^2-\alpha) + q_n (\alpha-1)$$
where $n=\sum_{i=2}^{N}\varepsilon_i R'_i,\; p_n=
-\sum_{i=2}^{N}\varepsilon_i G'_{i-2}$ and $q_n=
-\sum_{i=2}^{N}\varepsilon_i G'_{i-1}$. Using item (v) of
Proposition \ref{prop12}, we deduce  that if $x,y$ are the
coordinates of $\alpha^{2}$ in base $(\alpha^2- \alpha, \alpha-1)$,
then
 $1,\ x$ and $y$ are  $\mathbb{Q}$-Linearly independent.
Hence by Kronecker's Theorem, the set $\{n \alpha^{2}+p (\alpha^2-
\alpha) +q (\alpha-1),\ n \in \mathbb{N},\ p,\ q \in \mathbb{Z}\}$
is a dense set in $\mathbb{C}$. Using the fact that ${\mathcal R}_a$ is
the closure of the set $\{ \sum_{i=2}^{N}\varepsilon_i \alpha^i,\;
(\varepsilon_i)_{2 \leq i \leq N} \in E(R)\}$ and the same proof of
Proposition \ref{prop12}, we deduce that $\mathbb{C}=\bigcup_{v \in
\mathbb{Z} [\alpha-1]+\mathbb{Z}[\alpha^2- \alpha]} \mathcal{R} \cap
(\mathcal{R}+v)$ and if $(int(\mathcal{R})+v) \cap \mathcal{R}\neq
\emptyset$ where $v \in H= \mathbb{Z}
[\alpha-1]+\mathbb{Z}[\alpha^2- \alpha]$ then $v=0$.

On the other hand, the boundary of ${\mathcal R}_a$ is given by  $\partial {\mathcal R}_a = \bigcup_{v \in H -\{0\}}\mathcal{R}_a \cap (\mathcal{R}_a + v)$ .
Let $w \in S_a = \{0,\pm \alpha,\pm \alpha^2, \pm (\alpha- \alpha^2),\;\pm
(1+(a-1)\alpha^2), \pm (1+ (a-2)\alpha^2),\; \pm (1 -\alpha+
(a-1)\alpha^2),\;\pm (1-2 \alpha +a\alpha^2)\}$ be a state of the automaton ${\mathcal A}$. From Proposition \ref{egalite} and  relation (\ref{Ak}), we
have $\mathcal{R}_a \cap (\mathcal{R}_a + w / \alpha) \ne \emptyset$.
We remove the states $w=0, \pm \alpha,\pm \alpha^2, \pm
(1+(a-1)\alpha^2)$, because in this cases $w / \alpha \not \in G$.
Since $B$ is equal to the set of $w / \alpha$ such that $w \in \{ \pm (\alpha- \alpha^2),\;
 \pm (1+ (a-2)\alpha^2),\; \pm (1 -\alpha+
(a-1)\alpha^2), \pm (1-2 \alpha +a\alpha^2)\}= \pm (\alpha^3- \alpha^2),\;\pm (1-2 \alpha +a\alpha^2), \}$ , we obtain that
$ \bigcup_{v \in B} \mathcal{R}_a \cap (\mathcal{R}_a + v) \subset \partial{\mathcal R}_a .$

Now, let $z$ be an element of ${\mathcal R}_a$. Considering $\{\alpha^2- \alpha, \alpha-1\}$   instead of $\{1,\alpha\}$ and using exactly the same argument done
in the proof of item 2 of Proposition \ref{teofront2}, we deduce
that  there exist $(\varepsilon_{i})_{i \geq 2}$ and
$(\varepsilon_{i}^{'})_{i \geq l} \in  E(R),\ l<2$ such that
$z=\sum_{i=2}^{+\infty}
\varepsilon_{i}\alpha^{i}=\sum_{i=l}^{+\infty}
\varepsilon_{i}^{'}\alpha^{i}$ and $\varepsilon_{l}^{'} \neq 0$.
Hence
 $r=(\varepsilon_{l}^{'},0)
(\varepsilon_{l+1}^{'},0) \ldots
(\varepsilon_{1}^{'},0)(\varepsilon_{2}^{'}, \varepsilon_{2}) \ldots
,\; l <2$ is a path in the automaton ${\mathcal A}$ starting from the initial state. Therefore
$r$ satisfies one of the relations (\ref{eeqeq2002}), (\ref{eeqeq2003}), (\ref{eeqeq2004}), (\ref{eeqeq2005}).

If
$r = (1,0)(x_0,x_0 +1)(a-1,0)(x_1+a-2,x_1) (x_2,x_2+1) \ldots$, then $x_1=0$ and $z=(x_0+1)\alpha^{2}+\sum_{i=4}^{\infty} \varepsilon_{i}\alpha^{i} =
\alpha +
x_0\alpha^{2}+(a-1)\alpha^{3}+\sum_{i=4}^{\infty}\varepsilon_{i}^{'}\alpha^{i}= -(\alpha^2- \alpha)+ (x_0+1) \alpha^{2}+(a-1)\alpha^{3}+\sum_{i=4}^{\infty}\varepsilon_{i}^{'}\alpha^{i}$.
Therefore $z \in \mathcal{R}_a \cap (\mathcal{R}_a- (\alpha^2- \alpha))$.

If
$r = (1,0)(x_0,x_0)(a-1,0)(x_{1}+a-2,x_{1}) (x_{2},x_{2})(0,a-1) \ldots$ then, if $x_0 >0$, we deduce that
$z \in \mathcal{R}_a \cap (\mathcal{R}_a- (\alpha^2- \alpha))$.

If $x_0=0$ and $x_2>0$, we have $l \in \{-2, -1, 0, 1\}$. If $l=-2$, then $z \in \mathcal{R} \cap (\mathcal{R}- (\alpha^2- \alpha)$. If $l=-1$, $z \in \mathcal{R} \cap (\mathcal{R}- (\alpha- \alpha)$. If $l=0$, then   $z \in \mathcal{R} \cap (\mathcal{R}- (\alpha- 1)$. If $l=1$, then
$z \in \mathcal{R}_a \cap (\mathcal{R}_a- (\alpha^2- \alpha)$.

If $x_0= x_2=0$, then $l \in \{-3, -2, -1, 0, 1\}$. If l= -3, then
$z \in \mathcal{R}_a \cap (\mathcal{R}_a+ (\alpha- 1))$.

The cases $r=(1,0)(x,x)(a-1,0)(x_{1}+a-2,x_{1})
(x_{2},x_{2}+1) \ldots$  and $r=(1,0)(x,x)(x_1+a-2,x_1)
sss$ where  $s=(a-1,0)(0,a-1)(0,a-1)(a-1,0) $ are left to the reader.

Using the automaton $\mathcal{A}$, we deduce that $ \mathcal{R}_a
\cap (\mathcal{R}_a +  (1+ \alpha) u)= \{-1\},\;
 \mathcal{R}_a \cap (\mathcal{R}_a +  ( \alpha-1) u)= \{-\alpha\}$.
 Indeed, if $z\in \mathcal{R}_a
\cap (\mathcal{R}_a +  (1+ \alpha) u)=R_a\cap(R_a+\al^{-2}+\al)$,
then the path representing $z$ in the automaton is
$$(1,0)(0,0)(0,0)(1,0)(1,0)(0,1)ttt...,$$ where
$t=(1,0)(1,0)(0,1)(0,1)$. Hence, $z=-1$.\\
The case $\mathcal{R}_a \cap (\mathcal{R}_a + (\alpha-1)
u)=R_a\cap(R_a+\al^{-1})$ is left to the reader. \vspace{3em}
\hfill$\Box$


%

\section {Parametrization of the boundary of $\mathcal{R}_2$}

In this section, for simplicity, we consider the case $a=2$ and we give a complete description about the boundary of
$\mathcal{R}_2$ where

\[
\mathcal{R}_2=\left\{\sum_{i=2}^{\infty} \varepsilon_{i}\alpha^{i}
,\varepsilon_{i}\in\{0,1\}\ \forall i\geq 2,\ \varepsilon_{i}
\varepsilon_{i-1}\varepsilon_{i-2}\varepsilon_{i-3} <_{lex}
1101\right\}.
\]

We have seen (Proposition \ref{prop35}) that $\partial \mathcal{R}_2
= \bigcup_{v \in B} \mathcal{R}_2 \cap (\mathcal{R}_2+v)$ where
$B=\{\pm(\alpha^{-3}+\alpha^{-1})=\pm(\alpha-1),\ \pm \alpha^{-1},\
\pm (1+\alpha^{-2}),\ \pm (\alpha+\alpha^{-2})\}$. Since
$\mathcal{R}_2 \cap (\mathcal{R}_2 \pm v)$ is a point if $v=
\alpha^{-2}+\al$ or $v= \alpha^{-1}$. We will study the others four
regions $ \mathcal{R}_v = \mathcal{R}_2 \cap (\mathcal{R}_2 + v)$
where $v \in \{\pm (\alpha-1),  \pm ( 1+ \alpha^{-2})\}$, in
particular, we will prove the following result.

\begin{prop}
\label{lem51} Let $g$ and $h_i,\; i=0,1,2$ be the functions defined
by $g(z)= \al-1+\al z,\;  h_{0}(z)=\alpha-1+\alpha^{2}z,\
h_{1}(z)=-1+\alpha^{3}z$ and
$h_{2}(z)=\alpha^{2}+\alpha^{3}+\alpha^{4}z$ for all $z \in
\mathbb{C}$. Then we have the following properties:
\begin{enumerate}
\item $\mathcal{R}_{\al-1}=g(\mathcal{R}_{\al^{-2}+1})$.
\item $\mathcal{R}_{\alpha-1}= h_{0}(\mathcal{R}_{\alpha-1}) \cup h_{1}(\mathcal{R}_{\alpha-1}) \cup h_{2}(\mathcal{R}_{\alpha-1}).$
\item $h_{1}(\mathcal{R}_{\alpha-1}) \cap h_{2}(\mathcal{R}_{\alpha-1}) = \{-1-\alpha^{2}-\alpha^{4}\}=\{h_{2}(-\alpha-\alpha^{-1})\}.$
\item $h_{1}(\mathcal{R}_{\alpha-1}) \cap h_{0}(\mathcal{R}_{\alpha-1}) =
\{-1-\alpha^{3}\}=\{h_{1}(-1)\}.$
\item $h_{0}(\mathcal{R}_{\alpha-1}) \cap h_{2}(\mathcal{R}_{\alpha-1}) = \emptyset$.
\end{enumerate}
\end{prop}

\begin{obs}

Using  b), c), d)  and e) of the last Proposition, we will construct an explicit continuous and bijective application from $[0,1]$ to $\mathcal{R}_{\alpha-1}$.
Using this fact and a), we obtain an explicit homeomorphism between the circle and the boundary of ${\mathcal R}_2$.

\end{obs}

\begin{lem}\label{prop51} The
following properties are valid:

\begin{enumerate}

\item $\mathcal{R}_{\alpha^{-3}+\al^{-1}}\cap \mathcal{R}_{1+\alpha^{-2}}= \{-1\}.$
\item $\mathcal{R}_{\alpha^{-3}+\alpha^{-1}}\cap \mathcal{R}_{-1-\alpha^{-2}}= \{-\alpha-\alpha^{-1}\}.$
\item $\mathcal{R}_{\alpha^{-1}} \cap \mathcal{R}_{1+\alpha^{-2}} = \{-\alpha\}.$
\end{enumerate}
\end{lem}

\begin{obs}
 For the proof of Lemma \ref{prop51}, we will use the following relations:
\begin{equation}\label{rel}
\forall n \in \mathbb{Z},\;
\alpha^{n}=\alpha^{n-1}+\alpha^{n-2}+\alpha^{n-4},\;
\alpha^{n}+\alpha^{n-2}=2\alpha^{n-1}+\alpha^{n-3}.
\end{equation}
\end{obs}

\dem

\item[1)]
Let $z$ be an element of $\mathcal{R}_{\alpha^{-3}+\alpha^{-1}} \cap
\mathcal{R}_{1+\alpha^{-2}}$. The path  representing $z$ in the
automaton is $(1,0)(0,1)(1,0)(0,1)(0,0)(0,1)tttt...$ where
$t=(1,0)(1,0)(0,1)(0,1)$. Then
$$z=\alpha^{-3}+\alpha^{-1}+\sum_{1}^{\infty}(\alpha^{4i-1}+\alpha^{4i}).$$
By (\ref{rel}), we have
$z+\alpha=\al^{-3}+\al^{-1}$, then  $z=-1$.\\
\item[2)]
Let $z$ in $\mathcal{R}_{\alpha^{-3}+\alpha^{-1}} \cap
\mathcal{R}_{-1-\alpha^{-2}}$. As
$\al^{-3}+\al^{-1}+\al^{-2}+1=\al^{-2}+\al,$
 then $z+1+\al^{-2}\in\mathcal{R}_{\alpha+\alpha^{-2}}$. The path  representing $z+1+\al^{-2}$ in the
automaton is $(1,0)(0,0)(0,0)(1,0)tttt...$ where
$t=(0,1)(0,1)(1,0)(1,0)$. Then
$$z+1+\al^{-2}=\sum_{1}^{\infty}(\alpha^{4i-2}+\alpha^{4i-1}).$$
By (\ref{rel}), we get $z+1+\al^{-2}=-1$,then
$z=-2-\al^{-2}=-\al-\al^{-1}$.\\
\item[3)] If $z$ is an element of $\mathcal{R}_{\alpha^{-1}} \cap
\mathcal{R}_{1+\alpha^{-2}}$, then $z\in\mathcal{R}_{\alpha^{-1}}$
and the path representing $z$ in the automaton is
$(1,0)(0,0)(0,0)(1,0)tttt...$ where $t=(0,1)(0,1)(1,0)(1,0)$. Then
$$z=\sum_{1}^{\infty}(\alpha^{4i-1}+\alpha^{4i}).$$
By (\ref{rel}) we get $z+\al=0$, that is, $z=-\al$.
 \hfill$\Box$

\end{enumerate}


{\bf{Proof of Proposition \ref{lem51}}}

a)
Let $z\in\mathcal{R}_{\al^{-2}+1}$, then, using the automaton we
have $z=\al^{-2}+1+\sum_{i\geq2}a_i\al^{i}=\sum_{i\geq2}b_i\al^{i}$
where $a_3=0$. Then
$$g(z)=\al-1+\sum_{i\geq2}b_i\al^{i+1}\in\mathcal{R}+\alpha-1$$
and
$$g(z)= \al^{-1}+2\al-1+\sum_{i\geq2}a_i\al^{i+1}=\al^{2}+\sum_{i\geq2}a_i\al^{i+1}\in
\mathcal{R}.$$ We conclude that
$g(\mathcal{R}_{\al^{-2}+1})\subseteq \mathcal{R}_{\al-1}$.\\ Now
given $w\in\mathcal{R}_{\al-1}$, using the automaton, we have
$$w=\al-1+\sum_{i\geq3}a_i\al^{i}=\al^2+\sum_{i\geq3}b_i\al^{i},\
\mbox{where\ }b_4=0,$$ we have $w=g(z),z\in\mathcal{R}_{\al^{-2}+1}$
where
$z=\al^{-2}+1+\sum_{i\geq3}b_i\al^{i-1}=\sum_{i\geq3}a_i\al^{i-1}$.\\
\ \\ We conclude that $\mathcal{R}_{\al-1}\subseteq
g(\mathcal{R}_{\al^{-2}+1})$ and then $\mathcal{R}_{\al-1}=
g(\mathcal{R}_{\al^{-2}+1})$.

Let $z$ be an element of $\mathcal{R}_{\alpha-1}$ using the
automaton we conclude that\\ $z=\al-1+\displaystyle\sum_{i\geq3}
a_i\al^{i}$ and $z=\al^2+\displaystyle\sum_{i\geq3} b_i\al^{i}$.
\begin{enumerate}
\item[b)] Using (\ref{rel}), we have
\begin{equation}\label{h0}
\begin{array}{l}h_0(z)=\al^2+\displaystyle\sum_{i\geq3}a_i\al^{i+2}=\al-1+\al^4+\displaystyle\sum_{i\geq3} b_i\al^{i+2}.\\
\end{array}
\end{equation}
\begin{equation}\label{h1}
\begin{array}{l}h_1(z)=\al^2+\displaystyle\sum_{i\geq3}a_i\al^{i+3}=\al-1+\al^3+\al^4+\displaystyle\sum_{i\geq3}
b_i\al^{i+3}.\\
\end{array}
\end{equation}

\begin{equation}\label{h2}
\begin{array}{l}h_2(z)=h_2(\al-1+\displaystyle\sum_{i\geq3}
a_i\al^{i})=\al-1+\al^3+\al^4+\displaystyle\sum_{i\geq3}a_i\al^{i+4};\\
\\h_2(z)=h_2(\al^2+\displaystyle\sum_{i\geq3}
b_i\al^{i})=\al^2+\al^3+\al^6+\displaystyle\sum_{i\geq3}
b_i\al^{i+4}.
\end{array}
\end{equation}
Therefore $h_{i}(\mathcal{R}_{\alpha-1}) \subset
\mathcal{R}_{\alpha-1},\ \forall i \in \{0,1,2\}$ and then
$$h_{0}(\mathcal{R}_{\alpha-1}) \cup h_{1}(\mathcal{R}_{\alpha-1})
\cup h_{2}(\mathcal{R}_{\alpha-1})\subseteq \mathcal{R}_{\alpha-1}.$$

Now take $z\in \mathcal{R}_{\alpha-1}$.\\
If $(a_3,b_3)=(0,0)$ then $(a_4,b_4)=(1,0)$ and
$$z=\al-1+\al^4+\sum_{i\geq 5}a_i\al^{i}=\al^2+\sum_{i\geq5}b_i\al^{i}\ \mbox{with}\
a_7a_6a_51<_{lex}1101.$$ Using (\ref{h0}), we get $z=h_0(z_0)$ where
$z_0$ is the element of $\mathcal{R}_{\alpha-1}$ given by
$$z_0=\al-1+\sum_{i\geq5}b_i\al^{i-2}=\al^2+\sum_{i\geq5}a_i\al^{i-2}.$$

If $(a_3,b_3)=(1,0)$ then $(a_4,b_4)(a_5,b_5)=(1,0)(0,0)$ and
$$z=\al-1+\al^3+\al^4+\sum_{i\geq 6}a_i\al^{i}=\al^2+\sum_{i\geq6}b_i\al^{i}\ \mbox{with}\
a_8a_7a_61<_{lex}1101.$$ By (\ref{h1}) we get $z=h_1(z_0)$ where
$z_0$ is the element of $\mathcal{R}_{\alpha-1}$ given by
$$z_0=\al-1+\sum_{i\geq6}b_i\al^{i-3}=\al^2+\sum_{i\geq6}a_i\al^{i-3}.$$

If $(a_3,b_3)=(1,1)$ then
$(a_4,b_4)(a_5,b_5)(a_6,b_6)=(1,0)(0,0)(0,1)$ and
$$z=\al-1+\al^3+\al^4+\sum_{i\geq 7}a_i\al^{i}=\al^2+\al^3+\al^6+\sum_{i\geq7}b_i\al^{i}\ \mbox{ with }
b_9 b_8 b_71<_{lex}1101.$$
By (\ref{h2}), we have $z=h_2(z_0)$ where
$z_0$ is the element of $\mathcal{R}_{\alpha-1}$ given by
$$z_0=\al-1+\sum_{i\geq7}a_i\al^{i-4}=\al^2+\sum_{i\geq7}b_i\al^{i-4}.$$
Therefore $\mathcal{R}_{\al-1}\subseteq h_0(\mathcal{R}_{\al-1})\cup
h_1(\mathcal{R}_{\al-1})\cup h_2(\mathcal{R}_{\al-1}).$

c) Let  $z\in h_{1}(\mathcal{R}_{\alpha-1}) \cap
h_{2}(\mathcal{R}_{\alpha-1})$. Then there exist $z_1,\ z_{2} \in
\mathcal{R}_{\alpha-1}$ such that
$h_1(z_1)=-1+\alpha^{3}z_1=\alpha^{2}+
\alpha^{3}+\alpha^{4}z_2=h_2(z_2)$. Then $z_1=\al+\al z_2
=-1-\alpha^{-2}+ \alpha^{2}+\alpha z_2$. Since $z_{2} \in
\mathcal{R}_{\alpha-1}$, using the automaton,
$z_2=\al^2+\sum_{i=3}^{+\infty} b_i\al^{i},$ where $b_4=0.$
Therefore
$$z_1=-1-\alpha^{-2}+ \alpha^{2}+\alpha z_2=-1-\alpha^{-2}+\al^2+
\alpha^{3}+\sum_{i\geq3} b_i\alpha^{i+1}.$$ Hence
$z_1\in\mathcal{R}_{\alpha-1}\cap
\mathcal{R}_{-1-\al^{-2}}=\{-\al^{-1}-\al\}$ and then
$$z=h_1(z_1)=-1-\al^2-\al^4.$$

\item[d)] Let  $z\in h_{0}(\mathcal{R}_{\alpha-1}) \cap h_{1}(\mathcal{R}_{\alpha-1})$.
Then there exist $z_0,\ z_1 \in \mathcal{R}_{\alpha-1}$ such that
$h_0(z_0)=\al-1+\alpha^{2}z_0=-1+\alpha^{3}z_1=h_1(z_1)$ and then we
have $z_1=\al^{-2}+\al^{-1}z_0$. Since
$z_0\in\mathcal{R}_{\alpha-1}$, then
$z_0=\al^2+\sum_{i\geq3}b_i\al^{i}$ and
$$z_1=\al^{-2}+\al^{-1}z_0=\al^{-2}+\al^{-1}(\al^2+\sum_{i\geq3}b_i\al^{i})=\al^{-2}+\al+\sum_{i\geq3}b_i\al^{i-1}\in\mathcal{R}_{\al^{-2}+\al}.$$
Therefore
$z_1\in\mathcal{R}_{\alpha^{-2}+\al}\cap\mathcal{R}_{\al-1}=\{-1\}$.
Then
$$z=h_1(z_1)=h_1(-1)=-1-\al^{3}.$$

e) Left to the reader, can be done by the same manner that the
others items.

\hfill$\Box$
\end{enumerate}
\begin{lem}
\label{formu65} Consider $\ h_{2}$ as in Proposition \ref{lem51} and
$z\in\mathbb{C}$. Then  $\displaystyle \lim _{n \longrightarrow
\infty} h_{2}^{n}(z)=-1.$
\end{lem}
\dem
Since $h_2(z)=\al^2+\al^3+\al^4 z$, we can prove by induction that
$h_2^n(z)=\displaystyle\sum_{i=0}^{n-1}(\al^{4i+2}+\al^{4i+3})+\al^{4n}z$
for all integer $n \geq 1$. Hence
 $$\displaystyle \lim _{n \longrightarrow \infty}
h_{2}^{n}(z)=\sum_{i=0}^{\infty}(\al^{4i+2}+\al^{4i+3})=\frac{\al^2+\al^3}{1-\al^4}=-1.$$
%

\hfill$\Box$
\subsubsection{Parametrization of $\mathcal{R}_{\alpha-1}$}
Here, we will give an explicit parametrization of $\mathcal{R}_{\alpha-1}$
and hence for the boundary $\partial\mathcal{R}$. Let $z$ be an
element of $\mathcal{R}_{\alpha-1}$. Using Proposition \ref{lem51}, there
exists a sequence $(z_{n})_{n \geq 1}$ in $\mathcal{R}_{\alpha-1}$,
such that
$$
z=h_{a_{1}} \circ h_{a_{2}} \circ \ldots \circ h_{a_{n}}
(z_{n}),\forall n\geq1.
$$
If $x$ is an element of $\mathcal{R}_{\alpha-1}$, the sequence
$y_n=h_{a_{1}} \circ h_{a_{2}} \circ \ldots \circ h_{a_{n}}(x)$
converges to $z$ because the functions $h_{i}, i=0,1,2$ are
contractions.

Taking $x_{0} \in \mathcal{R}_{\alpha-1}$, $t\in[0,1],\
t=\sum_{i=1}^{\infty} a_{i} 3^{-i},a_{i}\in \{0,1,2\}$ we can define
a function $f:[0,1] \longrightarrow \mathcal{R}_{\alpha-1}$ by
$f(t)=\displaystyle \lim_{n \longrightarrow \infty} h_{b_{1}} \circ
\ldots \circ h_{b_{n}} (x_{0})$ where $(b_{i})_{i\geq1} =
\psi((a_{i})_{i\geq1})$, and

 $$\begin{array}{ccc}\psi:
\{0,1,2\}^{\mathbb{N}}&\longrightarrow&\{0,1,2\}^{\mathbb{N}}\\
a_1a_2...& \longmapsto& b_1b_2...\end{array}$$ is defined as follows: Put
$b_{1}=a_{1}$.

For $k \geq 2$, we define $b_{k}$ as follows:


\begin{equation}
 \mbox{ If }   a_{k}=1  \mbox{ then }  b_{k}=1 \label{eq547} ;
\end{equation}

\begin{equation}
\begin{array}{cccc}
\mbox{If } a_{k} \neq 1 \mbox{ and } & a_{k-1}=2 & \mbox{then} &
b_{k}=a_{k}; \label{eq548}
\end{array}
\end{equation}
\begin{equation}
\begin{array}{cccc}
\mbox{If } a_{k} \neq 1
 \mbox{ and } & a_{k-1}=0 &
\mbox{then } b_k =2 - a_k. \label{eq549}
\end{array}
\end{equation}

\vspace{0.3cm}

If $a_{k} \neq 1$ and $a_{k-1}=1$, let $r=min \{1 \leq i \leq k-1,\
a_{i}=a_{i+1}=\ldots= a_{k-1}=1\}$:

If ($r > 1$ and $a_{r-1}=2$) or $r=1$ then
\begin{equation}
\left\{\begin{array}{ll}  b_{k} = a_{k}, & \mbox{if} \hspace{0.2cm}
(k-r) \hspace{0.3cm} \mbox{is even} \\  b_{k}=2-a_{k}, & \mbox{if}
\hspace{0.2cm} (k-r) \hspace{0.3cm} \mbox{is odd}
 \end{array} \right. \label{eq5410}
\end{equation}

where $(k-r)$ is the number of digits $1$ after a number $0$ or $2$.

\vspace{0.3cm}

If $r > 1$ and $a_{r-1}=0$ then
\begin{equation}
\left\{\begin{array}{ll}b_{k}=a_{k}, & \mbox{if} \hspace{0.2cm}
(k-r) \hspace{0.3cm} \mbox{is odd} \\ b_{k}=2-a_{k}, & \mbox{se}
\hspace{0.2cm} (k-r) \hspace{0.3cm} \mbox{is even}
 \end{array} \right. \label{eq5411}
\end{equation}


\begin{teo}\label{teofh}: The application $f:[0,1] \longrightarrow \mathcal{R}_{\alpha-1}$ is well
defined, bijective, continuous and $f(0)=-\alpha-\alpha^{-1},\
f(1)=-1$.
\end{teo}

For the proof,  we need the following classical Lemma.

\vspace{0.3cm}

\begin{lem}\label{L4}
 Let $t$, $t^{'}$ be elements in $[0,1]$, $t=\sum_{i=1}^{\infty} a_{i}
 3^{-i},\
t^{'}=\sum_{i=1}^{\infty} c_{i} 3^{-i}$ with
$a_{i},c_{i}\in\{0,1,2\}$ such that $a_{i}=c_{i}$ for $i<k$ and
$a_{k} < c_{k}$ for some $k \in \mathbb{N}^{\star}$. Then
\begin{enumerate}
\item If $\mid t-t^{'} \mid < 3^{-N}, N > k$ then
$c_{k}=a_{k}+1$ and $\ c_{i}=0,\ a_{i}=2,\ k+1\leq i\leq N.$

\item If $t=t^{'}$ then $c_{k}=a_{k}+1$ and $\ c_{i}=0,\ a_{i}=2\ \forall i\geq k+1.$
\end{enumerate}
\end{lem}

\textbf{Proof of Theorem \ref{teofh}:}
Let $t,t^{'}\in[0,1]$ such that
$t=\sum_{i=1}^{\infty} a_{i} 3^{-i}$, $t^{'}=\sum_{i=1}^{\infty}
a'_{i} 3^{-i}$ with $a_{i},a'_{i}\in\{0,1,2\}$, $a_{i}=a'_{i}$ for
$i<k$, $a_{k} < a'_{k}$ for some integer $k \in \mathbb{N}$. \\
Assume that $f(t)=\displaystyle \lim_{n \longrightarrow \infty}
g_{b_{1}} \circ \ldots \circ g_{b_{n}} (x_{0})$,
$f(t^{'})=\displaystyle \lim_{n \longrightarrow \infty} g_{b'_{1}}
\circ \ldots \circ g_{b'_{n}} (x_{0})$ where $(b_{i})= \psi(a_{i})$
and $(b'_{i})= \psi(a'_{i}).$

{\bf $f$ is well defined:} Suppose that $t=t^{'}$ and
$(a_{i})_{i \geq 1} \leq_{lex} (a_{i}^{'})_{i \geq 1}$. There exists an integer
$k\in\mathbb{N}$ such that
\[
\begin{array}{ccc}
t=a_{1} \ldots a_{k-1} c \overline{2} & \mbox{and} & t^{'}= a_{1} \ldots a_{k-1} (c+1) \overline{0}
\end{array}
\] where $c=0$ or $1$,\;  $\overline{0}=00000...$ and $\overline{2}=2222...$.

If $t=a_{1} \ldots a_{k-1} 0 \overline{2}$ and $t^{'}= a_{1} \ldots
a_{k-1} 1 \overline{0}$ then $\psi(t)=b_{1}\ldots b_{k-1} 2 0
\overline{2},\ \psi(t^{'})=b_{1}\ldots b_{k-1} 1 0 \overline{2}$ if
$a_{k-1}=0$.
If $a_{k-1}=1$, then  $\psi(t)=b_{1}\ldots b_{k-1} x 0 \overline{2},\
\psi(t^{'})=b_{1}\ldots b_{k-1} 1 x \overline{2}$, where $x=0\
\mbox{or}\ 2$.
If $a_{k-1}=2$, then $\psi(t)=b_{1}\ldots b_{k-1} 0 0
\overline{2},\ \psi(t^{'})=b_{1}\ldots b_{k-1} 1 \overline{2}$.
By Lemma \ref{formu65} and the fact that $h_2\circ
h_0(-1)=h_1\circ h_0(-1)$, $h_{x}\circ h_0(-1)=h_1\circ h_{x}(-1)$
for $x=0\ \mbox{or } 2$ and $h_0^2(-1)=h_1(-1)$ we deduce that
$\psi(t)=\psi(t')$.

If $t=a_{1} \ldots a_{k-1} 1 \overline{2}$ and $t^{'}= a_{1} \ldots
a_{k-1} 2 \overline{0}$ then $\psi(t)=b_{1}\ldots b_{k-1} 1
\overline{2},\ \psi(t^{'})=b_{1}\ldots b_{k-1} 0 0 \overline{2}$ if
$a_{k-1}=0$, $\psi(t)=b_{1}\ldots b_{k-1} 1 x \overline{2},\
\psi(t^{'})=b_{1}\ldots b_{k-1} x 0 \overline{2}$, where $x=0\
\mbox{or}\ 2$ if $a_{k-1}=0$ and $\psi(t)=b_{1}\ldots b_{k-1} 1
\overline{2},\ \psi(t^{'})=b_{1}\ldots b_{k-1} 20 \overline{2}$ if
$a_{k-1}=2$. As before, we deduce that $\psi(t)=\psi(t')$.

\vspace{0.3cm}

{\bf $f$ is injective:} We know that $a_{i}=a_{i}^{'},\ 1 \leq
i \leq k-1$. Then\\
$f(t)=f(t^{'}) \Longleftrightarrow
h_{b_k}(z)=h_{b^{'}_k}(z')\Longleftrightarrow
(b_k=0,z=-\alpha-\alpha^{-1},b'_{k}=1,z'=-1)\ \mbox{or}\
(b_k=1,b'_k=2,z=z'=-\al-\alpha^{-1}$. Hence, we have to consider the
following cases:
\begin{itemize}
\item \textbf{Case 1} $b_k=0,\; b'_{k}=1,\; z=-\alpha-\alpha^{-1},\; z'=-1$.\\
According to Proposition \ref{lem51}, we have
$(b_i)=b_1b_2...b_{k-1}00\overline{2}$ and
$(b'_i)=b_1b_2...b_{k-1}1\overline{2}$.\\ If $b'_k=1$ then, by
 (\ref{eq547}), $a'_k=1$. As $b_k=0$ we have to consider the
following sub cases:
\begin{itemize}
\item \textbf{Case 1.1}: $a_{k-1}=2$.\\ By  (\ref{eq548}),
$a_k=0$. We have $b_{i}=2$ for all $i \geq k+1$. Hence by
(\ref{eq548}), $a_{i}=2,\ \forall i \geq k+1$. We also have
$b_{j}^{'}=2$  for all $j \geq k+1$. By (\ref{eq549}),
$a_{j}^{'}=0,\ \forall j \geq k+1$. Hence, $(a_i)_{i \geq
1}=a_1a_2...a_{k-2}20\overline{2}$ and $(a'_i)_{i \geq
1}=a_1a_2...a_{k-2}21\overline{0}$. Hence $t=t'$.
\item \textbf{Case 1.2}: $a_{k-1}=0$ and $b_{k}=0$.\\ We have $b_{k+1}=0$ and by
(\ref{eq549}) $a_k=2$, $b_{k}=0$ and $a_{k+1}=2$. As $b_{i}=2,\
\forall i \geq k+2$ then by (\ref{eq549}) we have $a_{i}=0$ for all
$i \geq k+2$. Since $b_{j}^{'}=2,\ \forall j \geq k+1$ then by
(\ref{eq548}) we have $a_{j}^{'}=2,\ \forall j \geq k+1$. Hence,
$(a_i)_{i \geq 1}=a_1a_2...a_{k-2}02\overline{0}$ and $(a'_i)_{i
\geq 1}=a_1a_2...a_{k-2}01\overline{2}$. Thus $t=t'$.
\item \textbf{Case 1.3}: $a_{r-1}=0,a_r=...=a_{k-1}=1$ and $(k-r)$ is
even.\\
If $b_{k}=0$ then by (\ref{eq5411}) $a_k=2$. Using (\ref{eq549}) and
 (\ref{eq548}) we have $(a_i)_{i \geq
1}=a_1...a_{r-2}011...12\overline{0}$ and $(a'_i)_{i \geq
1}=a_1...a_{r-2}011...11\overline{2}$. Hence $t=t'$.

\noindent Using the same arguments we prove the following cases:
$a_{r-1}=0,a_r=...=a_{k-1}=1$ and $(k-r)$ odd,
$a_{r-1}=2,a_r=...=a_{k-1}=1$, $(k-r)$ even,
$a_{r-1}=2,a_r=...=a_{k-1}=1$  $(k-r)$ odd, $a_1=a_2=...=a_{k-1}=1$
and $(k-1)$ even and $a_1=a_2=...=a_{k-1}=1$ and $(k-1)$ odd.

\end{itemize}
\item \textbf{Case 2} $b_k=1,\ b'_k=2,\ z=z'=-\alpha-\alpha^{-1}.$\\
According to Proposition \ref{lem51} we have
$(b_i)=b_1b_2...b_{k-1}10\overline{2}$ and
$(b'_i)=b_1b_2...b_{k-1}20\overline{2}$.\\ As $b_k=1$ then  $a_k=1$.
As $b'_{k+1}=0$ using the same previous ideas we have to consider
the following cases:
\begin{itemize}
\item \textbf{Case 2.1}: $a_{k-1}=2$.\\
We have $(a_i)_{i \geq 1}=a_1a_2...a_{k-2}21\overline{2}$ and
$(a'_i)_{i \geq 1}=a_1a_2...a_{k-2}22\overline{0}$. Then
\ref{L4}, $t=t'$.
\item \textbf{Case 2.2}: $a_{k-1}=0$.\\
We have $(a_i)_{i \geq 1}=a_1a_2...a_{k-2}01\overline{0}$ and
$(a'_i)_{i \geq 1}=a_1a_2...a_{k-2}00\overline{2}$. Then
 $t=t'$.
\item \textbf{Case 2.3}: $a_{r-1}=0,a_r=...=a_{k-1}=1$ and $(k-r)$ even.\\
Here we have $(a_i)_{i \geq 1}=a_1...a_{r-2}011...11\overline{0}$
and $(a'_i)_{i \geq 1}=a_1...a_{r-2}011...10\overline{2}$. Then $t=t'$.

Using the same arguments we prove the following cases:
$a_{r-1}=0,a_r=...=a_{k-1}=1,\ (k-r)$ odd,
$a_{r-1}=2,a_r=...=a_{k-1}=1,\ (k-r)$ even,
$a_{r-1}=2,a_r=...=a_{k-1}=1,\ (k-r)$ odd, $a_1=a_2=...=a_{k-1}=1,\
(k-1)$ even and $a_1=a_2=...=a_{k-1}=1,\ (k-1)$ odd.

\end{itemize}
\end{itemize}

{\bf $f$ is continuous:} Let
$t=\sum_{i=1}^{\infty}a_{i}3^{-i}$,
$t^{'}=\sum_{i=1}^{\infty}a_{i}^{'}3^{-i}$ and suppose that $0 <
\mid t - t^{'} \mid < 3^{-N},\ N \in \mathbb{N},\; N > k.$ By Lemma
\ref{L4}, $a_{k}^{'}=a_{k}+1,\ a_{i}^{'}=0$ and $a_{i}=2$ for all
$i$ satisfying $k+1 \leq i \leq N.$ We have to consider the
following cases:
\begin{itemize}
\item \textbf{Case 1.1:} If $a_{k-1}=0,\ a_{k}=0,\ a_{k}^{'}=1$ then
we can write
$$t= (a_{1} \ldots a_{k-2} 0 0  \overline{2} ),\  t^{'} = (a_{1}
\ldots a_{k-2} 0 1 \overline{0})$$

and since $|h_i(z)-h_i(w)|\leq |\al|^2|z-w|$ then
\[
\mid f(t)- f(t^{'}) \mid = \mid f(t)-f(t^{'}) \mid =\mid h_{2}(z_{1})-h_{1}(z_{1}^{'}) \mid
. \mid \alpha \mid ^{2(k-1)}.
\]

Where $z_{1}=h_{0}h_{2}^{N-k-1}(y_{1})$ and
$z_{1}^{'}=h_{0}h_{2}^{N-k-1}(y_{1}^{'}), y_1, y_{1}^{'} \in \mathbb{C}$
\[
\mid f(t)-f(t^{'}) \mid =\mid h_{2}(z_{1})-h_{1}(z_{1}^{'}) \mid
. \mid \alpha \mid ^{2(k-1)}
\]

As
$h_{2}(-\alpha-\alpha^{-1})=h_{1}(-\alpha-\alpha^{-1})=-(1+\alpha^{2}+\alpha^{4})$ then
\[
\begin{array}{ccc}
\mid f(t)-f(t^{'}) & \leq & \mid
h_{2}(z_{1})-h_{2}(-\alpha-\alpha^{-1})\mid + \mid
h_{1}(-\alpha-\alpha^{-1}) - h_{1}(z_{1}^{'}) \mid \times \mid
\alpha \mid ^{2(k-1)}\\ & \leq & (1+\mid \alpha\mid) \mid
\alpha\mid^{2k+1} diam ({\mathcal R}_{\alpha -1}),
\end{array}
\] where $diam ({\mathcal R}_{\alpha -1})$ is the diameter of ${\mathcal R}_{\alpha
-1}$.

\item \textbf{Case 1.2:} If $a_{k-1}=2,\ a_{k}=0,\ a_{k}^{'}=1$ then
we can write
$$
t= (a_{1} \ldots a_{k-2} 2 0  \overline{2} ),\ t^{'} = (a_{1} \ldots
a_{k-2} 2 1 \overline{0} ). $$ Then
\[
\mid f(t)-f(t^{'}) \mid =\mid h_{2}(z_{2})-h_{1}(z_{2}^{'}) \mid .
 \mid \alpha \mid ^{2(k-1)},
\]

where $z_{2}=h_{0}h_{2}^{N-k-1}(y_{2})$ and
$z_{2}^{'}=h_{0}h_{2}^{N-k-1}(y_{2}^{'}), \; y_2, y´_{2} \in \mathbb{C}$.
Hence
\[
\begin{array}{ccc}
\mid f(t)-f(t^{'})  & \leq & (1+\mid \alpha \mid) \mid \alpha
\mid^{2k+1} diam({\mathcal R}_{\alpha -1}).
\end{array}
\]

\item \textbf{Case 1.3:} If $a_{i-1}=0,\ a_{i}=\ldots=a_{k-1}=1,\
a_{k}=0,\ a_{k}^{'}=1$ and $(k-i-1)$ is even . This case can be done by the some as before and is left to the reader.

\end{itemize}

\hfill $\Box$



\section{ Hausdorff Dimension of $\partial(\mathcal{R}_2)$}

Since $\partial(\mathcal{R}_2)$ is the union of 4 curves that are
images of $\mathcal{R}_{\alpha-1}$ by a affine application, we have
that $dim_{H}\mathcal{R}_{\alpha-1}=dim_{H}\partial\mathcal{R}$. By
Proposition \ref{lem51}, the set $\mathcal{R}_{\alpha-1}=
\cup_{i=0}^{2}h_{i}(\mathcal{R}_{\alpha-1})$ is   invariant by the
affine maps  $h_{i}$. An upper bound of Hausdorff dimension of this
class of sets is given by Theorem

\vspace{0.3cm}

\begin{teo}\label{teo431}\cite{falconer} Let $A$ a set of $\mathbb{C}$ such that
$A=\cup_{i=0}^{n} \varphi_{i}(A)$ is compact and invariant for affine applications $\varphi_{i}$ with
coefficients $r_{i}$ (i.e. , $\forall x,\ y \in \mathbb{C},\ \mid
\varphi_{i}(x)-\varphi_{i}(y) \mid = r_{i} \mid x-y \mid$), then $dim_{H}(A)
\leq s$, where $s$ is the unique real number that verifies $\sum_{i=0}^{n}
r_{i}^{s}=1$.
\end{teo}

\begin{obs} When the $\varphi_{i}(A)$ intercept in points is known that $dim_{H}(A)=s$ (see \cite{falconer}).
\end{obs}

\vspace{0.3cm}

By Proposition \ref{lem51}, we deduce that
 $dim_{H}(\partial(\mathcal{R}_2))=s,$ where $s$
verifies
\[
\mid \alpha \mid^{2s} + \mid \alpha \mid^{3s} + \mid \alpha
\mid^{4s} =1.
\]

Therefore $dim_{H}(\partial(\mathcal{R}_2))=\frac{\log \rho}{\log
\mid \alpha \mid} \approx 1.359337357,$ where $\rho$ is the root is
the maximum real root of the polynomial $X^{4}+X^{3}+X^{2}-1=0$.

\begin{obs}
Using the automaton $\mathcal{A}$, we can also parametrize
$\mathcal{G}_2$ and prove that it is homeomorphic to a topological disk. We can also show that
$dim_{H}(\partial(\mathcal{G}_2))=dim_{H}(\partial(\mathcal{R}_2))$.
All these results can be extended to all $\mathcal{R}_a$ and $\mathcal{G}_a,\; a \geq 3$.
\end{obs}

\newpage

\begin{figure}[!htb]
\begin{minipage}[b]{0.55\linewidth}
\includegraphics[width=\linewidth]{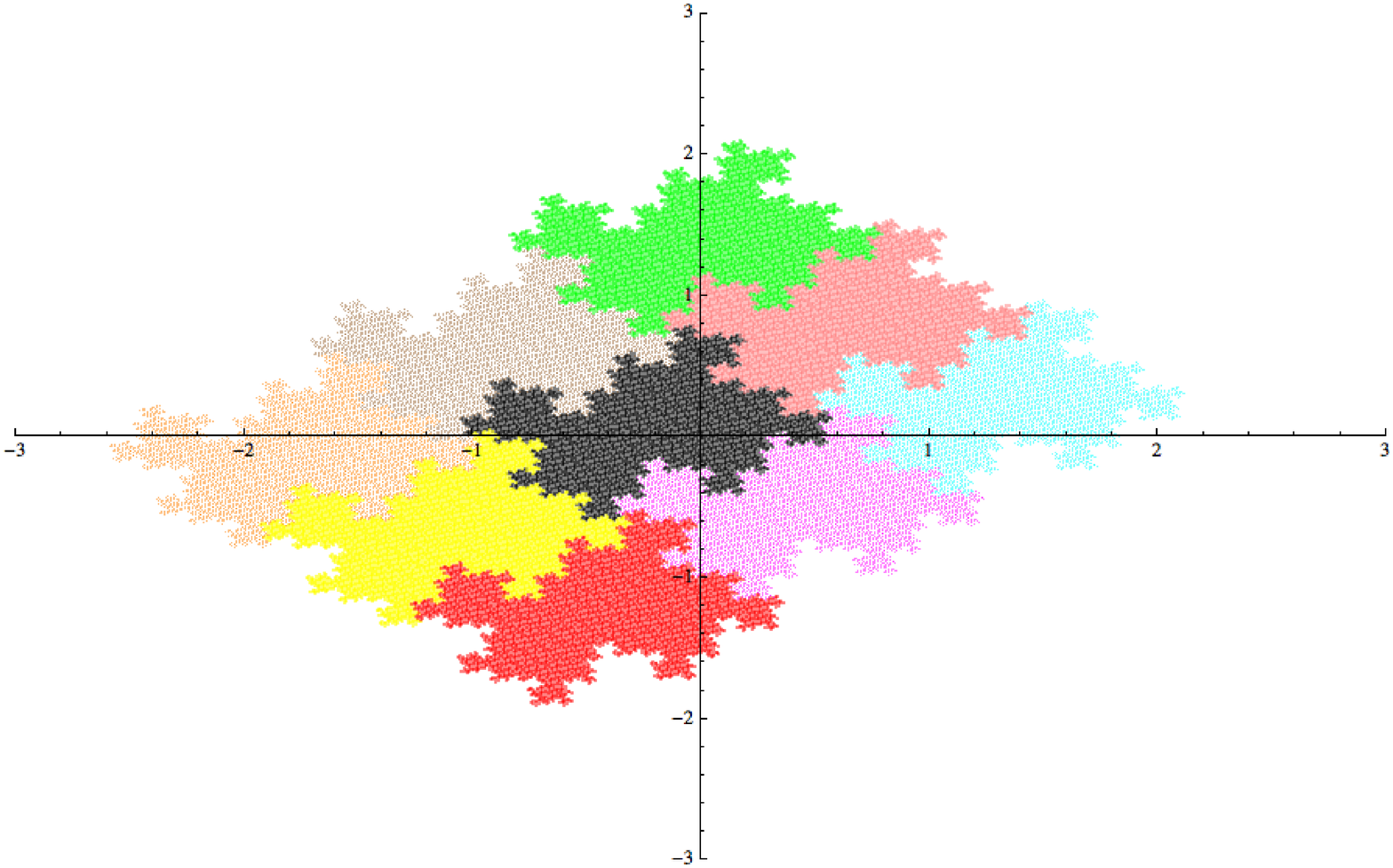}
\caption{Fractal $\mathcal{R}_{a}$} \label{fig:patu}
\end{minipage} \hfill
\begin{minipage}[b]{0.50\linewidth}
\includegraphics[width=\linewidth]{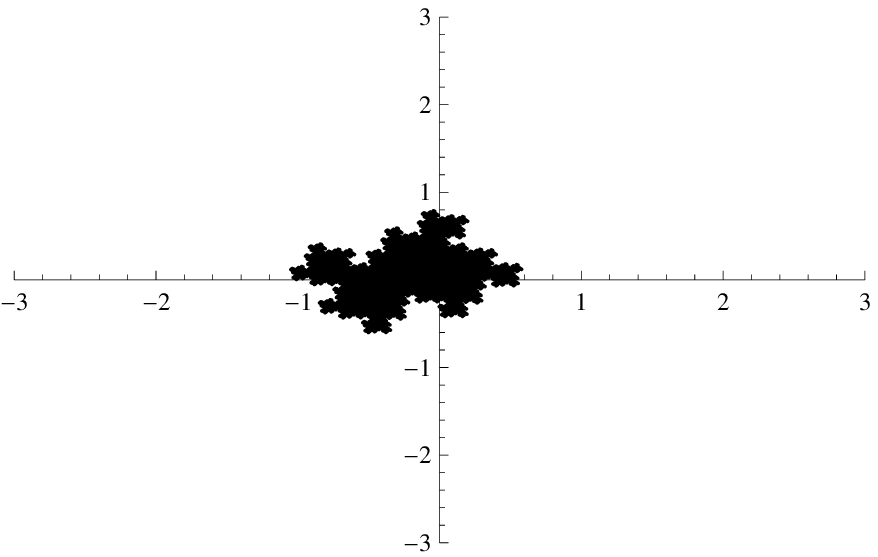}
 \caption{A tile of Fractal
$\mathcal{R}_{a}$} \label{fig:catole}
\end{minipage}
\end{figure}

\begin{figure}[h!]
    \begin{center}
        \includegraphics[scale=0.20]{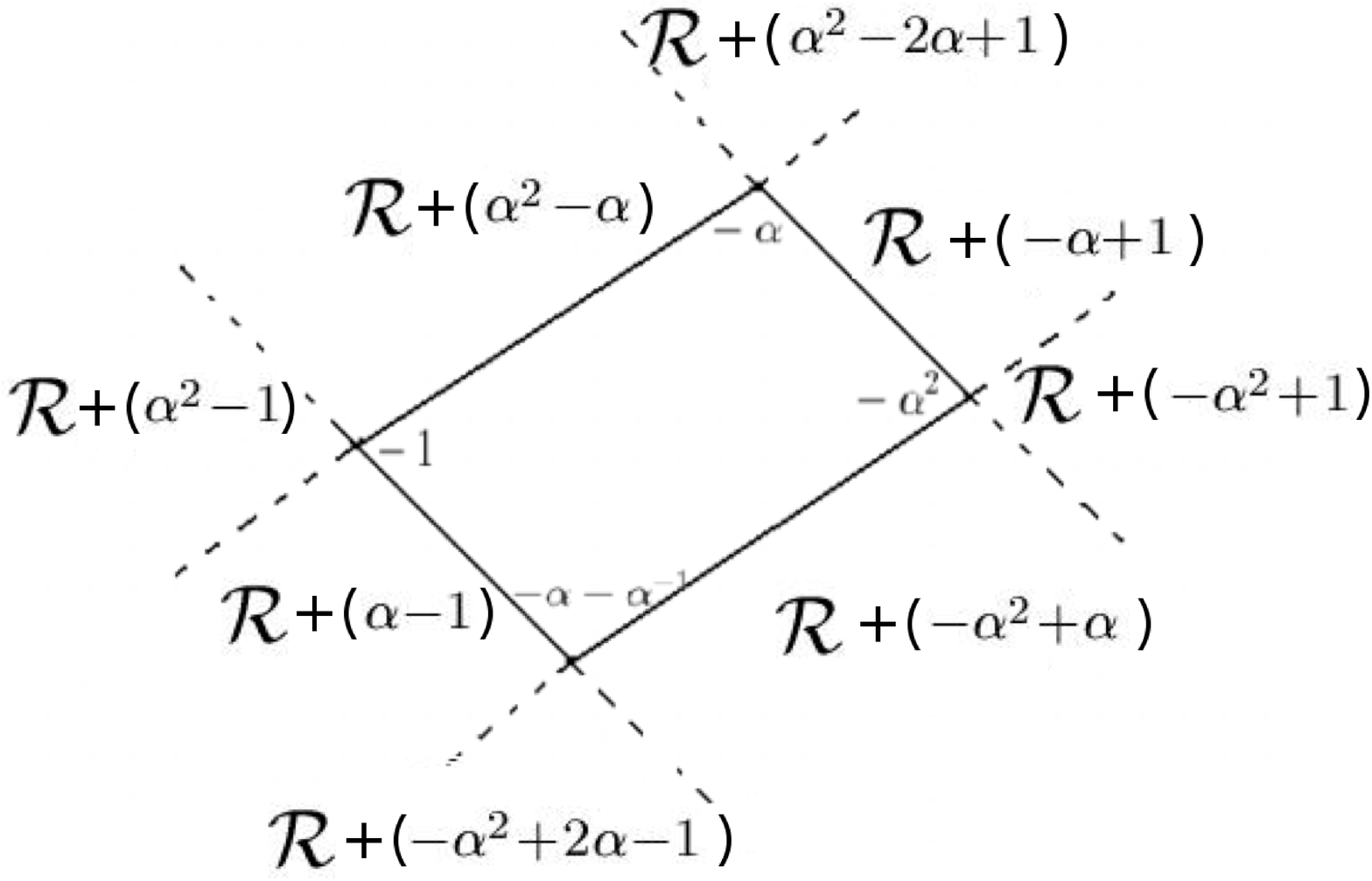}
        \caption{Boundary's Fractal $\mathcal{R}_{a}$}
    \end{center}
\end{figure}

\begin{figure}[h!]
\begin{center}
\includegraphics[scale=0.20]{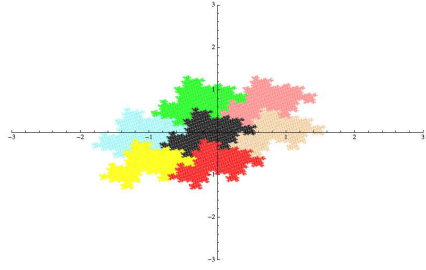}
\caption{Fractal $\mathcal{G}_{a}$}
\end{center}
\end{figure}

\begin{figure}[h!]
\begin{center}
\includegraphics[scale=0.20]{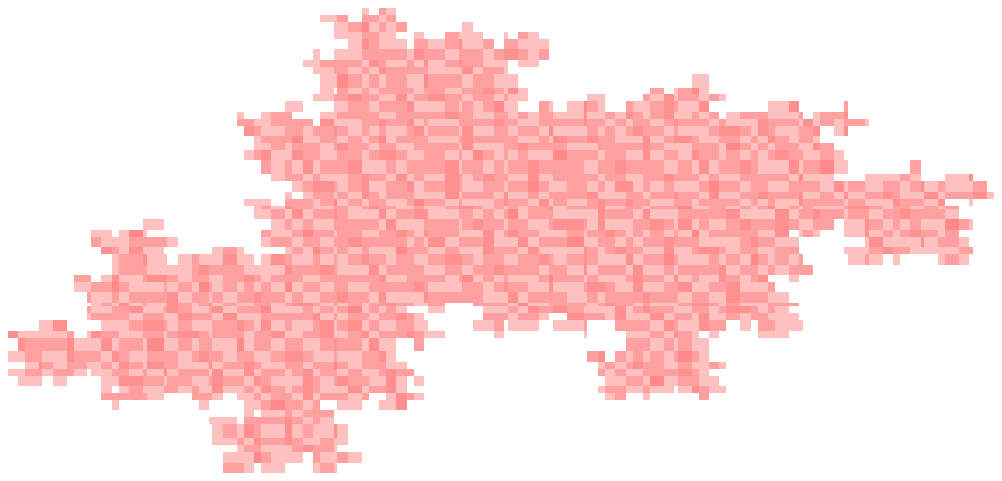}
 \caption{A tile of Fractal $\mathcal{G}_{a}$}
\end{center}
\end{figure}

\begin{figure}[h!]
\begin{center}
\includegraphics[scale=0.25]{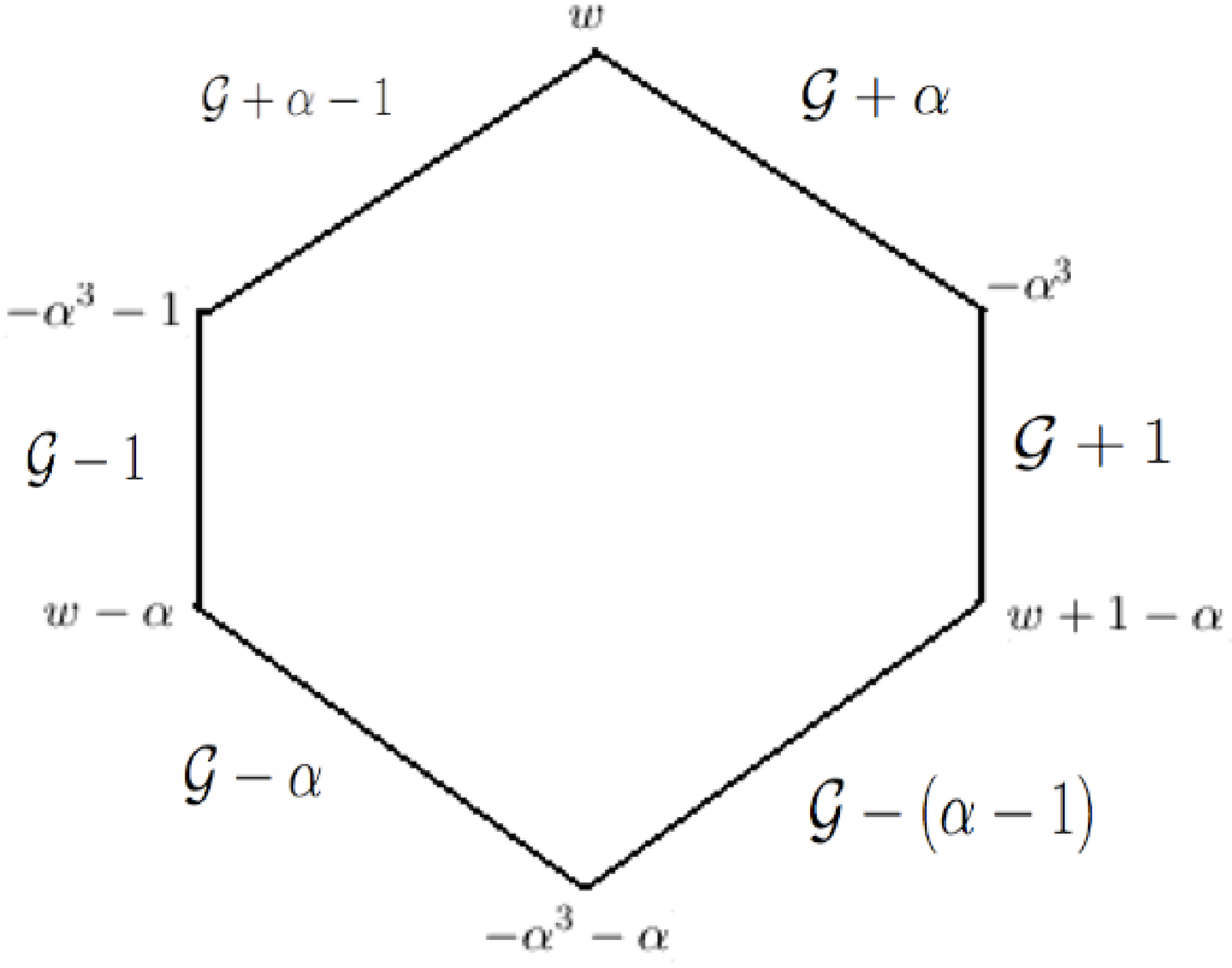}
 \caption{Boundary's Fractal $\mathcal{G}_{a}$}
\end{center}
\end{figure}


\end{document}